\documentclass[article]{ij4uq_arxiv}
\usepackage{eurosym}
\usepackage[pdftex]{thumbpdf}
\usepackage[pdftex,
pdfstartview=FitH,bookmarks=true,bookmarksnumbered=true,bookmarksopen=true,hypertexnames=false,breaklinks=true,
colorlinks=true,  linkcolor=blue,anchorcolor=blue,
citecolor=blue,filecolor=blue,  menucolor=blue]{hyperref}
\usepackage{amsmath}
\usepackage{amsfonts}
\usepackage{amssymb}
\usepackage{graphicx}%
\providecommand{\U}[1]{\protect\rule{.1in}{.1in}}
\providecommand{\U}[1]{\protect\rule{.1in}{.1in}}

\fancypagestyle{plain}{
\fancyhf{}
\fancyhead[R]{\small{\it}}
\fancyfoot[R]{\vspace*{10pt}\small\bf\thepage}
\fancyfoot[L]{\fottitle}
}
\begin{document}
\volume{Volume 1, Number 1, \myyear\today}
\title{Stochastic {G}alerkin finite element method
for nonlinear elasticity and application to reinforced concrete members}
\titlehead{Stochastic {G}alerkin {FEM} for nonlinear elasticity}
\authorhead{M. S. Ghavami \& B. Soused\'{\i}k \& H. Dabbagh \& M. Ahmadnasab}
\author[1]{Mohammad S. Ghavami}
\author[2]{Bed\v{r}ich Soused\'{\i}k}
\corrauthor[1]{Hooshang Dabbagh}
\author[3]{Morad Ahmadnasab}
\corremail{\href{mailto:h.dabbagh@uok.ac.ir}{h.dabbagh@uok.ac.ir}}
\address[1]{Department of Civil Engineering \\
University of Kurdistan \\
Sanandaj 66177-15175, Iran}
\address[2]{Department of Mathematics and Statistics \\
University of Maryland, Baltimore County \\
Baltimore, MD 21250, USA}
\address[3]{Department of Mathematics \\
University of Kurdistan \\
Sanandaj 66177-15175, Iran}
\dataO{Version date: \mydate\today}
\dataF{}
\abstract{
We develop a stochastic Galerkin finite element method for nonlinear elasticity
and apply it to reinforced concrete members with random material properties.
The strategy is based on the modified Newton-Raphson method, which consists of an incremental loading process
and a linearization scheme applied at each load increment.
We consider that the material properties are given by a stochastic expansion in the so-called generalized
polynomial chaos (gPC) framework.
We search the gPC expansion of the displacement,
which is then used to update the gPC expansions of the stress, strain and internal forces.
The proposed method is applied to a reinforced concrete beam with uncertain initial concrete modulus of elasticity \textcolor{black}{and a shear wall 
with uncertain maximum compressive stress of concrete},
and the results are compared to those of stochastic collocation and Monte Carlo methods.
Since the systems of equations obtained in the linearization scheme
using the stochastic Galerkin method are very large, 
and there are typically many load increments, 
we also studied iterative solution using preconditioned conjugate gradients.
The efficiency of the proposed method is illustrated by a set of numerical experiments.
}
\maketitle

\section{Introduction}

\label{sec:introduction}Reinforced concrete is one of the most popular
materials in construction. Therefore predicting the behavior of reinforced
concrete members and structures is very important. However, development of the
analytical models for reinforced concrete is complicated by several factors.
Reinforced concrete is composed of concrete and steel, which are two materials
with different nonlinear properties. In particular, concrete exhibits a
nonlinear softening behavior in compression and a linear brittle behavior in
tension, and its combination with steel makes the behavior of reinforced
concrete members even more complex. The computational analysis is commonly
performed using the finite element method. To this end, considerable number of
constitutive models was developed to characterize the nonlinear behavior and
stress strain relationship of concrete and its reinforcement. In traditional
approaches, the physical characteristics of reinforced concrete are considered
to be known and the problem is deterministic. However, in practice the
structural properties of these materials typically show variability, which
usually result from the manufacturing process, natural variability in
microstructure and possibly also from aging. Additional factors, such as
uncertainty in external loading,
may also contribute to the uncertainty in the predicted response.
Propagation and quantification of
uncertainty using numerical simulations for risk and reliability analysis as
well as design of reinforced concrete members and structures is therefore a
critical task of engineering research.

The most widely applied technique for approximating quantities of interest for
problems with random inputs is Monte Carlo method~\cite{Fishman1996montecarlo,
caflisch1998monte}. It is simple to implement, but the main weakness is
relatively slow convergence. Another important class is given by perturbation
methods~\cite{kleiber1992stochastic, kaminski2013stochastic}, which are
however limited to problems with small variability of uncertainty. An
alternative that gained a significant attention in the last two decades is the
stochastic finite element
method~\cite{ghanem1991stochastic,stefanou2009stochastic}. The assumption is
that a parametric uncertainty is described in terms of polynomials of random
variables using the so-called generalized polynomial chaos (gPC)
framework~\cite{xiu2002wiener} and one searches for gPC expansions of
solutions. There are two main approaches: \emph{stochastic collocation method}
(SC), which is based on sampling that translates the problem into a set of
uncoupled deterministic problems cf,
e.g.,~\cite{Babuska-2010-SCM,berveiller2006stochastic,gunzburger2014stochastic}%
, and \emph{stochastic Galerkin method} (SG), which by means of the Galerkin
projection couples the physical and probabilistic degrees of freedom into a
single large system of equations cf. also,
e.g.,~\cite{Babuska-2004-GFE,gunzburger2014stochastic,Matthies-2005-GML}.
Both approaches have been used in various applications in structural
engineering~\cite{sudret2015polynomial}, see also~\cite{Giraldi-2014-TBN}. Due
to the large size of the linear systems arising in the stochastic Galerkin
method, iterative solution may be the preferred choice since use of direct
solvers may be prohibitive. When the finite element discretization of the
underlying deterministic problem leads to a symmetric, positive definite
matrix, the global stochastic Galerkin matrix is typically also positive
definite~\cite{ghosh2014probabilistic}. Then the corresponding linear system
can be solved using the conjugate gradient method, and its preconditioning is
often a vital component enabling to speed up convergence and reduce
computational time. The most simple, yet often surprisingly effective is the
mean-based, block-diagonal preconditioner proposed by Pellissetti and
Ghanem~\cite{pellissetti2000iterative} and analyzed by Powell and
Elman~\cite{powell2009block}. Other preconditioners use either more terms from
the gPC expansion of the coefficient matrix or imitate the structure of the
stochastic Galerkin
matrix~\cite{ullmann2010kronecker,ullmann2012efficient,sousedik2014hierarchical,sousedik2014truncated}%
, or use advanced solvers as their components such as domain
decomposition~\cite{sarkar2009domain,ghosh2009feti,Subber-2014-DDM,Subber-2014-SPS}
or multigrid~\cite{rosseel2010iterative,Brezina-2014-SAA,Osborn-2015-MSS}, see
also~\cite{keese2005hierarchical,ernst2010stochastic,coulier2017inverse}. One
of the most recent advancements includes the truncation
preconditioners~\cite{sousedik2014truncated,kubinova2020block,Bespalov-2020-TPS}%
, which we also use in our numerical experiments. We note that a similar study
to ours was recently presented by~\cite{Lacour-2020-SFE}, but here we focus on
a more complex nonlinear model that includes a crack development, and we also
study efficient solution of the linear systems obtained by use of the
stochastic Galerkin method.

In this paper, we develop the stochastic Galerkin finite element method for
propagation of uncertainty in a nonlinear elasticity model and apply to
reinforced concrete members with random material properties. In particular, we
consider uncertainty in the modulus of elasticity, which (nonlinearly) depends
on the strain (deformation). We describe a linearization scheme based on the
modified Newton-Raphson method formulated in the context of the stochastic
Galerkin framework. The method consists of an incremental loading process and
a linearization scheme applied at each load increment. The tangent stiffness
matrix is updated after each load increment, and each increment is further
subdivided into a number of steps. In each such step, the structural response
is translated into the vector of internal forces, which then modifies the
right-hand sides of the linearized problems. The whole process continues until
the full load is applied and the displacement increment vector is smaller than
a given tolerance. We consider that the material properties (stiffness) are
given by a stochastic expansion, in the so-called generalized polynomial chaos
(gPC) framework. We search the gPC expansion of the displacement, which is
then used to update the gPC expansions of the stress, strain and internal
forces between each step of a load increment. Stochastic tangent stiffness
matrices are updated at the first step of each load increment. The 
method is applied to a reinforced concrete beam with uncertain initial concrete modulus of elasticity 
\textcolor{black}{and a shear wall with uncertain maximum compressive stress of concrete} 
under deterministic concentrated loading and the results
are compared to that of Monte Carlo and stochastic collocation methods. Next,
since the size of the linear systems in the stochastic Galerkin method is
usually large, we also study iterative solution using preconditioned conjugate
gradients. 

The paper is organized as follows. In Section~\ref{sec:dfem} we recall the
general implementation of deterministic nonlinear FEM, in
Section~\ref{sec:sfem:galerkin} we formulate the stochastic Galerkin finite
element method that accounts for a nonlinearity in the material properties, in
Section~\ref{sec:sfem:sampling} we discuss the sampling methods (Monte Carlo
and stochastic collocation), in Section~\ref{sec:sfem:example} we present the
application of the proposed procedure to a simply supported beam 
\textcolor{black}{and in Section~\ref{sec:sfem:example2} to a shear wall}
with random material properties, In
section~\ref{sec:precond} we study iterative solution of the linear systems
arising in the stochastic Galerkin method, and finally in
Section~\ref{sec:conclusion} we summarize and conclude our work.

\section{ Deterministic nonlinear model}

\label{sec:dfem} We recall the deterministic nonlinear finite element
formulation for reinforced concrete members based on~\cite{Bazant-1982-SRF}.
We consider
that the
response is mainly caused by nonlinearity in stress-strain relations, due to
cracking of concrete and due to yielding of steel reinforcement. In general,
structure
of concrete is very complex because it is a composite made up of hydrated
cement, sand, and coarse aggregates. In addition, it contains numerous flaws
and micro-cracks, and the rapid propagation of micro-cracks under applied
loads contributes to the nonlinear properties of concrete as well. Numerical
simulation of reinforced concrete members requires realistic stress-strain
relations of the plain concrete as an input. Such formulations are thus
nonlinear and rely on a number of experimentally determined material constants.

For the concrete, we will use the representation by the nonlinear elasticity
model based on the limiting tensile strain failure criterion that covers the
pre-peak and post-peak regimes from~\cite{kupfer1969behavior}. This model is
expressed in terms of tangent stiffness formulation. It is assumed that
stresses in the principal directions can be calculated independently of each
other, based on uniaxial stress-strain relationships. The biaxial effect was
assumed to be due to the interaction of the two principal directions through
the Poisson's ratio effect. Therefore, the concept of \emph{equivalent
uniaxial strain} should be introduced in order to separate the Poisson effect
from the cumulative strain~\cite{darwin1977nonlinear}. In addition, equivalent
uniaxial strain can make a contribution to keep track of the degradation of
stiffness and strength of plain concrete and to allow actual biaxial
stress-strain curves to be derived from \emph{uniaxial} curves. The uniaxial
curves selected for compressive loading in this study are based on an equation
suggested by the FIB model code~\cite{code2010fib} as
\begin{equation}
\frac{\sigma_{c}}{f_{cm}}=-\frac{k\eta-\eta^{2}}{1+(k-2)\eta}\ \quad
\text{for}\ |\epsilon_{c}|<|\epsilon_{c,lim}|, \label{equation.conccomp}%
\end{equation}
where $\sigma_{c}$ is the compressive stress in MPa, $f_{cm}=\acute{f}_{c}+8$
is the actual compressive strength of concrete at an age of 28 days in MPa
(here $\acute{f}_{c}$ is the specific compressive strength of concrete in
MPa), $k=E_{ci}/E_{c1}$ is the plasticity number (here $E_{ci}$ is the tangent
modulus of elasticity in MPa, and $E_{c1}$ is the secant modulus from the
origin to the peak compressive stress), and $\eta=\epsilon_{c}/\epsilon_{c1}$
(here $\epsilon_{c}$ is the compressive strain and $\epsilon_{c1}$ is the
strain at the maximum compressive stress). Then, a complete stress-strain
curve of concrete in uniaxial tension can be treated based on an equation
suggested by FIB model code~\cite{code2010fib} as
\begin{equation}
\sigma_{ct}=%
\begin{cases}
E_{ci}\,\epsilon_{ct} & \text{for}\ \sigma_{ct}\leq0.9f_{ctm},\\
f_{ctm}\left(  1-0.1\frac{0.00015-\epsilon_{ct}}{0.00015-0.9f_{ctm}/E_{ci}%
}\right)  & \text{for}\ 0.9f_{ctm}<\sigma_{ct}\leq f_{ctm}%
\end{cases}
\label{equation.concten}%
\end{equation}
where $\sigma_{ct}$ is the tensile stress in MPa, $E_{ci}$ is the tangent
modulus of elasticity in MPa, $\epsilon_{ct}$ is the tensile strain, and
$f_{ctm}$ is the tensile strength in MPa.

In contrast with concrete, the mechanical properties of steel reinforcement
are well known. The reinforcing steel bars are assumed to be only transmitting
axial compressive or tensile forces. Thus, an uniaxial stress-strain
relationship is sufficient. Under monotonic loading, it is generally assumed
that the steel behavior is identical in tension and compression. In this
study, a bilinear stress-strain relation suggested by FIB model
code~\cite{code2010fib} is used as
\begin{equation}
\sigma_{s}=%
\begin{cases}
E_{s}\,\epsilon_{s} & \text{for}\ \epsilon_{s}<\epsilon_{sy},\\
f_{y}+E_{sh}\,(\epsilon_{s}-\epsilon_{sy}) & \text{for}\ \epsilon_{sy}%
\leq\epsilon_{s}\leq\epsilon_{su}%
\end{cases}
\label{equation.steel}%
\end{equation}
where $\sigma_{s}$ is the stress, $E_{s}$ is the initial modulus of
elasticity, $\epsilon_{s}$ is the strain, $\epsilon_{sy}$ is the yielding
strain, $f_{y}$ is the yielding stress, $E_{sh}$ is the modulus of strain
hardening and $\epsilon_{su}$ is the ultimate strain.

The phenomenon of concrete cracking is extremely important in the behavior of
reinforced concrete structures. The maximum stress and strain theories are
frequently used to determine whether tensile cracking has occurred in the
concrete~\cite{chen1982plasticity}. If the maximum principal stress or strain
in a point of the structure reaches the uniaxial tensile strength or tensile
strain limit, cracking is assumed to form perpendicular to the direction of
the maximum tensile stress or strain. The stress in that direction is
subsequently reduced to zero. The limiting tensile stress value which causes
the crack is not a well-defined quantity. For specimens cast from the same
concrete, the flexural tensile strength determined from a modulus of rupture
test is higher than the tensile strength of a split cylinder, which is in turn
higher than the tensile strength obtained from a direct tension test.
Furthermore, for each type of test there is significant scatter in the
results. For concrete structures subjected to rapid loading, the maximum
strain criterion is more realistic, since uniaxial dynamic tensile tests
indicate that an almost constant failure strain is observed irrespective of
the strain rate or loading rate~\cite{hatano1960dynamical}. The limiting
tensile strain criterion has been employed with success to represent the
tensile cracking of concrete under static loading~\cite{chen1980end}. 
\textcolor{black}{In this study, we use a smeared crack approach proposed by~\cite{rashid1968ultimate}.
It is based on the limiting tensile strain criterion, and it allows the concrete to crack in one or two directions. 
In fact, it offers complete generality in possible crack direction. 
This representation is also popular since it allows for automatic generation of cracks without redefinition of the
finite element topology. Instead of redefining the finite element mesh and nodal connectivities, 
cracking is accounted for by modifying the material properties within elements. 
After the crack has occurred in an element, the concrete becomes an orthotropic material 
with one of the material axes being oriented along the direction of cracking,
and the elasticity modulus in the direction perpendicular to the cracking plane is reduced to zero.
Next, a} reduced shear modulus is assumed on the cracked plane to account
for aggregate interlocking. Shear transmission due to aggregate interlock can
be simply accounted for in the smeared crack model by the introduction of a
reduced value of concrete shear modulus. In this study, we used a proposed
model by~\cite{al1979nonlinear}. The tension stiffening effect of cracked
concrete is also incorporated into this model by including a descending branch
in the stress-strain curve of concrete under tension. In this study, we use a
linear descending branch in the stress-strain curve of concrete under tension
for considering tension stiffening effect, and other nonlinear effects such as
crushing of concrete in compression and yielding or strain hardening of steel
reinforcement are also taken into account, as suggested by FIB model
code~\cite{code2010fib}.

We combine the constitutive equations with the compatibility and equilibrium
equations, compute the weak form of the equilibrium equations and replace the
displacement field by its finite element approximation. This yields a discrete
nonlinear equilibrium system, which may be formally written as
\begin{equation}
\mathcal{K}(u)=f, \label{equation.nonlinearequequ}%
\end{equation}
where $\mathcal{K}(u)$ is the nonlinear operator that depends on the model
parameters but also on the displacement vector $u$, and $f$ is the vector of
external load. We now formulate the modified Newton-Raphson method to solve
system~(\ref{equation.nonlinearequequ}). The method is based on linearization
of the nonlinear operator and incremental application of the external load.
The tangent stiffness matrix is updated after each load increment, and the
system is assumed linear during a load increment. Each load increment is
further subdivided into several steps, during which only the right-hand side
is updated. The system of linearized equations at step$~i$ of the $n$th load
increment can be written as
\begin{equation}
K^{n}\Delta u^{n,i}=f^{n}-g\left(  u^{n,i-1}\right)  ,
\label{equation.linearequequ}%
\end{equation}
where $K^{n}$ is the corresponding tangent stiffness matrix of size~$n_{x}%
\times n_{x}$, $\Delta u^{n,i}$ is the increment of displacement vector,
$f^{n}$ is the external load vector and $g\left(  u^{n,i-1}\right)  $ is the
internal force vector at $(i-1)$th step of the $n$th load increment. The
matrix $K^{n}$ is obtained by assembly of finite element matrices. The matrix
of the $e$th element at the $n$th load increment is computed as
\begin{equation}
K_{e}^{n}=\int_{v_{e}}B_{e}^{T}D_{e}^{n}B_{e}\,dv_{e}, \label{equation.kne1}%
\end{equation}
where $B_{e}$ is the element shape function derivative matrix and $D_{e}^{n}$
is the element stress-strain matrix computed as
\begin{equation}
D_{e}^{n}=D_{c,e}^{n}+D_{s,e}^{n}, \label{equation.dne}%
\end{equation}
where $D_{c,e}^{n}$ and $D_{s,e}^{n}$ are concrete and steel rebar
stress-strain matrices, respectively. For the first load increment, the
concrete stress-strain matrix is computed as
\begin{equation}
D_{c,e}^{1}=\left[
\begin{array}
[c]{ccc}%
\frac{E_{ci}}{1-\nu_{c}^{2}} & \frac{E_{ci}\nu_{c}}{1-\nu_{c}^{2}} & 0\\
\frac{E_{ci}\nu_{c}}{1-\nu_{c}^{2}} & \frac{E_{ci}}{1-\nu_{c}^{2}} & 0\\
0 & 0 & \frac{E_{ci}}{2(1+\nu_{c})}%
\end{array}
\right]  , \label{equation.d1ce}%
\end{equation}
where $E_{ci}$ is the initial concrete modulus of elasticity and $\nu_{c}$ is
the Poisson coefficient of concrete. For the first load increment, the steel
rebar stress-strain matrix is computed\ as
\begin{equation}
D_{s,e}^{1}=E_{si}\left[
\begin{array}
[c]{ccc}%
\rho_{x} & 0 & 0\\
0 & \rho_{y} & 0\\
0 & 0 & 0
\end{array}
\right]  , \label{equation.d1se}%
\end{equation}
where $E_{si}$ is the initial steel rebar modulus of elasticity and $\rho_{x}$
and $\rho_{y}$ are the steel rebar percentages in $x$ and $y$ directions,
respectively. For other load increments, $D_{c,e}^{n}$ and $D_{s,e}^{n}$ are
computed using the material constitutive models and the corresponding stress
level~\cite{Bazant-1982-SRF}. Applying numerical integration to
eq.~(\ref{equation.kne1}), we get
\begin{equation}
K_{e}^{n}=\sum_{q_{1}=1}^{m_{1}}\sum_{q_{2}=1}^{m_{2}}\omega_{q_{1}}%
\omega_{q_{2}}\left[  B_{e}(\zeta_{q_{1}},\eta_{q_{2}})\right]  ^{T}D_{e}%
^{n}\,B_{e}(\zeta_{q_{1}},\eta_{q_{2}})\det J(\zeta_{q_{1}},\eta_{q_{2}})\,t,
\label{equation.kne2}%
\end{equation}
where $m_{1}$, $m_{2}$ are the numbers of quadrature points $\zeta_{q_{1}},$
$\eta_{q_{2}}$ with weights $\omega_{q1}$, $\omega_{q2}$, respectively, $\det
J(\zeta_{q1},\eta_{q2})$ is the determinant of the Jacobian matrix, and $t$ is
the thickness of the element. The internal force vector$~g_{e}(u^{n,i-1})$ is
computed for each element in eq.~(\ref{equation.linearequequ})\ as
\begin{equation}
g_{e}(u^{n,i-1})=\int_{v_{e}}B_{e}^{T}\sigma_{e}^{n,i-1}\,dv_{e},
\label{equation.pui-11}%
\end{equation}
where $\sigma_{e}^{n,i-1}$ is the stress vector for the $e$th element.
Applying numerical integration to eq.~(\ref{equation.pui-11}), we get
\begin{equation}
g_{e}(u^{n,i-1})=\sum_{q_{1}=1}^{m_{1}}\sum_{q_{2}=1}^{m_{2}}\omega_{q_{1}%
}\omega_{q_{2}}[B_{e}(\zeta_{q_{1}},\eta_{q_{2}})]^{T}\sigma_{e}^{n,i-1}\det
J(\zeta_{q_{1}},\eta_{q_{2}})\,t. \label{equation.pui-12}%
\end{equation}
The stress vector $\sigma_{e}^{n,i-1}$ is computed as
\begin{equation}
\sigma_{e}^{n,i-1}=\sigma_{e}^{n,i-2}+\Delta\sigma_{e}^{n,i-1},
\label{equation.sigmai-1}%
\end{equation}
where $\Delta\sigma_{e}^{n,i-1}$ is the stress vector increment
\begin{equation}
\Delta\sigma_{e}^{n,i-1}=D_{e}^{n}\,\Delta\epsilon_{e}^{n,i-1},
\label{equation.deltasigmai-1}%
\end{equation}
and $\Delta\epsilon_{e}^{n,i-1}$ is the strain vector increment
\begin{equation}
\Delta\epsilon_{e}^{n,i-1}=B_{e}\,\Delta u_{e}^{n,i-1}.
\label{equation.deltaepsiloni-1}%
\end{equation}
Above, $\Delta u_{e}^{n,i-1}$ is the displacement increment at step~$i-1$ of
the $n$th load increment for the $e$th element. At the first step of the first
load increment, the deformation vector $u^{1,0}=0$, and so the internal force
vector $g(u^{n,i-1})=0$. After eq.~(\ref{equation.linearequequ}) is solved for
the displacement increment$~\Delta u^{n,i}$, a new approximate solution is
obtained as
\begin{equation}
u^{n,i}=u^{n,i-1}+\Delta u^{n,i}. \label{equation.ui}%
\end{equation}
This process continues until the displacement increment vector is smaller than
a given tolerance.

\section{Stochastic finite element method}

\label{sec:sfem}

\subsection{The stochastic Galerkin method}

\label{sec:sfem:galerkin} We assume that the uncertainty is induced in the
model by a vector$~\xi$ of independent, identically distributed (i.i.d.)
random variables$~\xi_{i}$, $i=1,\dots,m_{\xi}$. Specifically, we let the
uncertainty enter the model through some of the parameters, more details are
given in Section~\ref{sec:sfem:example}. Then, the nonlinear
system~(\ref{equation.nonlinearequequ}) becomes stochastic, and in this study
we use the stochastic Galerkin finite element framework to extend the
Newton-Raphson method from Section~\ref{sec:dfem} to this case.
This\ framework entails use of a space spanned by a set of multivariate
polynomials $\left\{  \psi_{\ell}\left(  \xi\right)  \right\}  _{\ell
=1}^{n_{\xi}}$, which is known in the literature as a generalized polynomial
chaos (gPC) basis~\cite{Xiu-2002-WAP,xiu2010numerical}. The polynomials are
orthonormal with respect to the density function associated with the
distribution of$~\xi$, and we will in particular assume that
\begin{equation}
\psi_{1}=1,\quad\text{and}\quad\mathbb{E}\left[  \psi_{k}\psi_{\ell}\right]
=\langle\psi_{k}\psi_{\ell}\rangle=\delta_{k\ell}, \label{equation.psi}%
\end{equation}
where $\mathbb{E}$ is the mathematical expectation, and$~\delta_{k\ell}$
denotes the Kronecker delta function.

Suppose we are given a stochastic\ expansion of the tangent stiffness matrix
at the $n$th load increment as
\begin{equation}
K(\xi)^{n}=\sum_{\ell=1}^{n_{K}}K_{\ell}^{n}\,\psi_{\ell}(\xi
),\label{equation.stiffdisc}%
\end{equation}
where $K_{\ell}^{n}$ is for each $\ell$ a (deterministic) matrix of size
$n_{x}\times n_{x}$. Let us further suppose that the external load vector at
the$~n$th load increment$~f(\xi)^{n}$ and the internal force vector at the
$(i-1)$th step of $n$th load increment$~g(\xi)^{n,i-1}$ are also given as
stochastic expansions
\begin{align}
f(\xi)^{n} &  =\sum_{m=1}^{n_{\xi}}f_{m}^{n}\,\psi_{m}(\xi
),\label{equation.loaddisc}\\
g(\xi)^{n,i-1} &  =\sum_{m=1}^{n_{\xi}}g_{m}^{n,i-1}\,\psi_{m}(\xi
),\label{equation.forcedisc}%
\end{align}
where both $f_{m}^{n}$ and $g_{m}^{n,i-1}$ are vectors of length~$n_{x}$ for
all$~m$.\ Forming the stochastic counterpart of
system~(\ref{equation.linearequequ}), we will search for the increment of the
displacement vector at the $i$th step of the $n$th load increment in the form
\begin{equation}
\Delta u(\xi)^{n,i}=\sum_{k=1}^{n_{\xi}}\Delta u_{k}^{n,i}\,\psi_{k}%
(\xi).\label{equation.dispdisc}%
\end{equation}
Specifically, substituting expansions~(\ref{equation.stiffdisc}%
)--(\ref{equation.dispdisc}) into equation~(\ref{equation.linearequequ}) and
performing a stochastic Galerkin projection, i.e., by orthogonalizing the
residual to the gPC basis $\{\psi_{m}\}_{m=1}^{n_{\xi}}$, yields a system of
linear deterministic equations
\begin{equation}
\mathbf{K}^{n}\Delta\mathbf{u}^{n,i}=\mathbf{f}^{n}-\mathbf{g}^{n,i-1}%
,\qquad\mathbf{K}^{n}\in\mathbb{R}^{n_{x}n_{\xi}\times n_{x}n_{\xi}}%
\quad\Delta\mathbf{u}^{n,i},\mathbf{f}^{n},\mathbf{g}^{n,i-1}\in
\mathbb{R}^{n_{x}n_{\xi}}.\label{equation.linearequequ2}%
\end{equation}
Denoting $c_{\ell km}=\langle\psi_{\ell}\psi_{k}\psi_{m}\rangle$, the
stochastic Galerkin matrix$~\mathbf{K}^{n}$ can be written as
\begin{equation}
\mathbf{K}^{n}=\left[
\begin{array}
[c]{ccc}%
K_{(1,1)}^{n} & \cdots & K_{({n}_{\xi},1)}^{n}\\
\vdots & \ddots & \vdots\\
K_{(1,{n}_{\xi})}^{n} & \cdots & K_{({n}_{\xi},{n}_{\xi})}^{n}%
\end{array}
\right]  ,\qquad\text{with }{K}_{(k,m)}^{n}=\sum_{\ell=1}^{n_{K}}c_{\ell
km}\,{K}_{\ell}^{n},\label{equation.stochstiff}%
\end{equation}
and the vectors in~(\ref{equation.linearequequ2}) are concatenations
of$~n_{\xi}$ subvectors of size$~n_{x}$, cf.~(\ref{equation.loaddisc}%
)--(\ref{equation.dispdisc}), as
\begin{equation}
\Delta\mathbf{u}^{n,i}=\left[
\begin{array}
[c]{c}%
\Delta u_{1}^{n,i}\\
\vdots\\
\Delta u_{n_{\xi}}^{n,i}%
\end{array}
\right]  ,\qquad\mathbf{f}^{n}=\left[
\begin{array}
[c]{c}%
f_{1}^{n}\\
\vdots\\
f_{n_{\xi}}^{n}%
\end{array}
\right]  ,\quad\mathbf{g}^{n,i-1}=\left[
\begin{array}
[c]{c}%
g_{1}^{n,i-1}\\
\vdots\\
g_{n_{\xi}}^{n,i-1}%
\end{array}
\right]  .\label{equation.stochdisp}%
\end{equation}
For the first step of the first load increment $\mathbf{g}^{1,0}=0$ since
$u^{1,0}(\xi)=0$. After solving equation~(\ref{equation.linearequequ2}), the
approximate solution is updated as
\begin{equation}
\mathbf{u}^{n,i}=\mathbf{u}^{n,i-1}+\Delta\mathbf{u}^{n,i}%
.\label{equation.ui2}%
\end{equation}
For the next steps, the gPC expansion of internal load vector$~g^{n,i-1}$ is
computed using numerical integration as
\begin{equation}
g_{m}^{n,i-1}=\left\langle g(\xi)^{n,i-1}\,\psi_{m}\right\rangle =\sum
_{q=1}^{n_{q}}g(\xi^{(q)})^{n,i-1}\,\psi_{m}(\xi^{(q)})\,\omega^{(q)}%
,\label{equation.puim}%
\end{equation}
where $\xi^{(q)}$ are quadrature points, $\omega^{(q)}$ are quadrature
weights, $n_{q}$ is the number of quadrature points and $g(\xi^{(q)})^{n,i-1}$
is a realization of expansion~(\ref{equation.forcedisc}) at $\xi^{(q)}%
$,\ which is obtained by assembly of finite element internal load vectors,
cf.~(\ref{equation.pui-12}),
\begin{equation}
g(\xi^{(q)})_{e}^{n,i-1}=\sum_{q_{1}=1}^{m_{1}}\sum_{q_{2}=1}^{m_{2}}%
\omega_{q_{1}}\omega_{q_{2}}[B(\zeta_{q_{1}},\eta_{q_{2}})]^{T}\sigma
(\xi^{(q)})_{e}^{n,i-1}\det J(\zeta_{q_{1}},\eta_{q_{2}}%
)\,t.\label{equation.puiq2}%
\end{equation}
The realizations of the element internal load vectors in~(\ref{equation.puiq2}%
) depend on the realizations of the stress vector~ $\sigma(\xi^{(q)}%
)_{e}^{n,i-1}$, which in turn depend on the realizations of the stress and
strain increments and are computed from the realizations of the displacement
increment~$\Delta u(\xi^{(q)})_{e}^{n,i-1}$ in the same way as the
calculations in~(\ref{equation.sigmai-1})--(\ref{equation.deltaepsiloni-1}).

The terms in expansion of the tangent stiffness
matrix~(\ref{equation.stiffdisc}) are also obtained using numerical
integration as
\begin{equation}
K_{j}^{n}=\left\langle K(\xi)^{n}\,\psi_{j}(\xi)\right\rangle =\sum
_{q=1}^{n_{q}}K(\xi^{(q)})^{n}\,\psi_{j}(\xi^{(q)})\,\omega^{(q)},
\label{equation.knj1}%
\end{equation}
where the realizations of $K(\xi^{(q)})^{n}$ are computed using realizations
of the model parameters, with more details given in
Section~\ref{sec:sfem:example}, and using the realizations of the displacement
vector~$u(\xi^{(q)})^{n,0}$ from its gPC\ expansion. This process is repeated
until the $\Delta u^{n,i}$ is smaller than a given tolerance. Statistics of
the solution, such as mean, standard deviation and probability density
function, can be estimated in post-processing phase by sampling the gPC expansions.

\subsection{Sampling methods}

\label{sec:sfem:sampling} In numerical experiments described in
Section~\ref{sec:sfem:example}, we compare results from stochastic Galerkin
method to those obtained using two sampling methods: Monte Carlo and
stochastic collocation. The sampling entails solving of a number of mutually
independent deterministic problems. In the Monte Carlo method, the sample
points $\xi^{(q)}$, $q=1,...,n_{MC}$, are generated randomly following the
distribution of the underlying random variables, and moments of the solution
are computed by ensemble averaging. For stochastic collocation, which is used
here in the form of so-called non-intrusive (or pseudospectral) stochastic
Galerkin method, the sample points $\xi^{(q)}$, $q=1,...,n_{q}$, are given as
a set of predetermined collocation points. This method derives from a sparse
grid for performing interpolation or quadrature using a small number of points
in multidimensional space~\cite{gerstner1998numerical, novak1996high}. There
are two procedure to implement stochastic collocation method to obtain the
coefficient in the equation~(\ref{equation.dispdisc}), the first way is
constructing a Lagrange interpolating polynomials and the second way is
performing a discrete projection using the so-called pseudospectral
approach~\cite{xiu2010numerical}. In this study we use the second approach for
a direct comparison with the stochastic Galerkin method. The coefficients of
the displacement vector at the $n$th load increment are computed using
numerical quadrature as
\begin{equation}
u_{k}^{n}=\left\langle u(\xi)^{n}\,\psi_{k}(\xi)\right\rangle =\sum
_{q=1}^{n_{q}}u(\xi^{(q)})^{n}\,\psi_{k}(\xi^{(q)})\,\omega^{(q)}.
\label{equation.deltauik}%
\end{equation}
This, in particular, means that instead of forming and
solving~(\ref{equation.linearequequ2}), the increments\ of the displacement
vector $\Delta u(\xi^{(q)})^{n,i}$ are found for each sample point$~\xi^{(q)}$
independently from
\begin{equation}
K(\xi^{(q)})^{n,i}\,\Delta u(\xi^{(q)})^{n,i}=f^{n}-g(\xi^{(q)})^{n,i-1},
\label{equation.ktnq}%
\end{equation}
where $K(\xi^{(q)})^{n,i}$ and $g(\xi^{(q)})^{n,i-1}$ are realization of the
stochastic tangent stiffness matrix and internal forces vector at quadrature
points, respectively, and $f^{n}$ is the external load vector.

\subsection{Example~1: Simply supported reinforced concrete beam}

\label{sec:sfem:example} 
\begin{figure}[b]
\centering
\includegraphics[width=12cm,height=7cm]{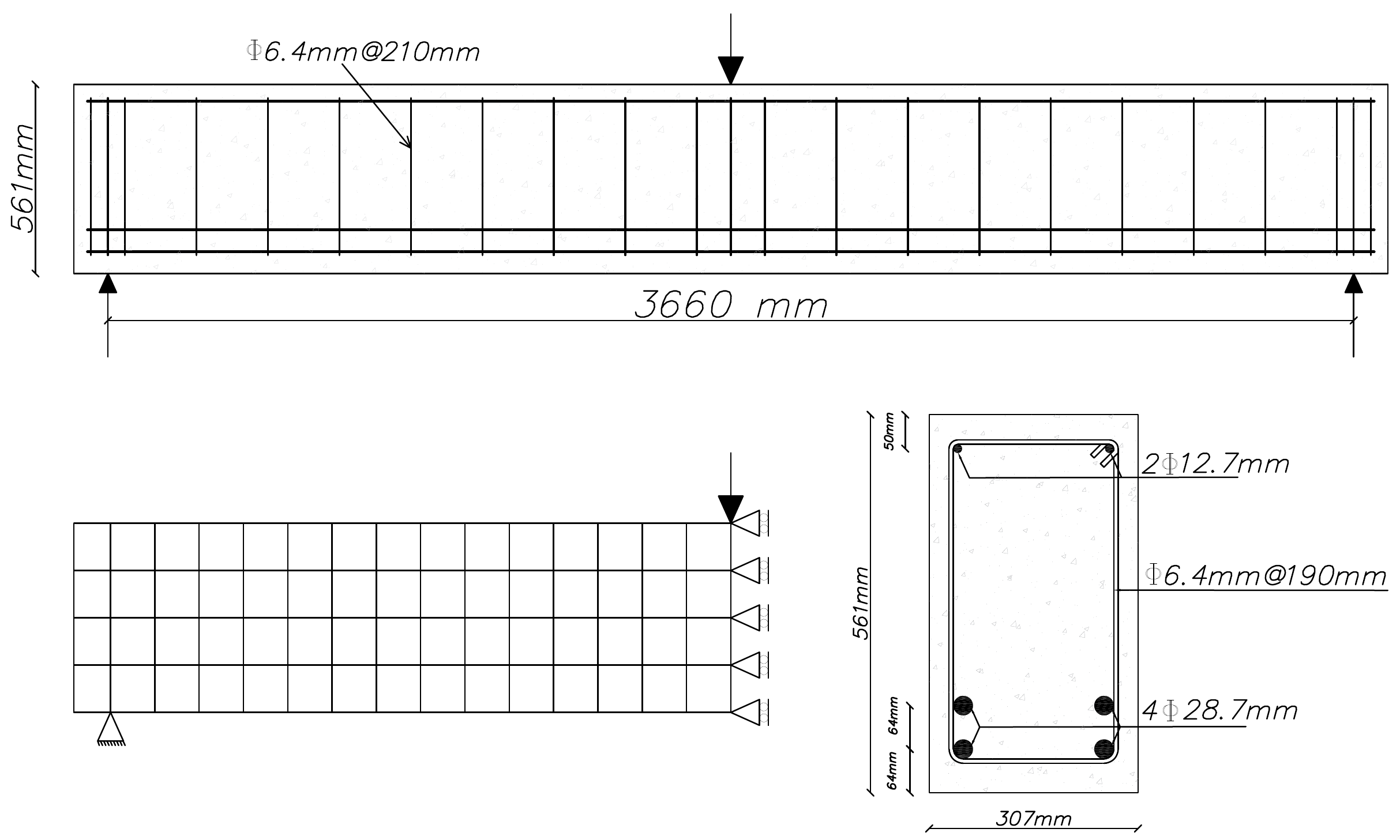}
\caption{Technical drawing of the beam. }%
\label{figure.beama1}%
\end{figure}

\begin{figure}[b]
\centering
\includegraphics[width=12cm]{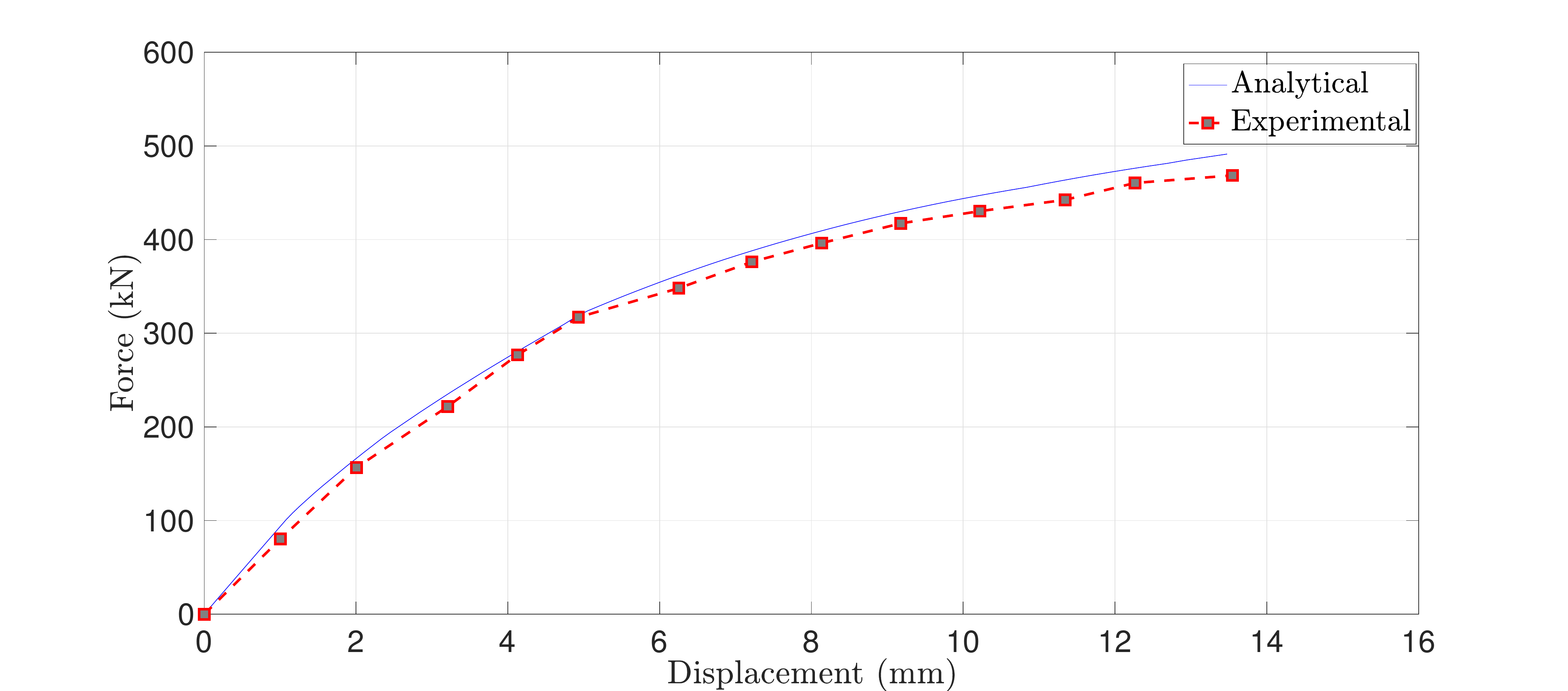} \caption{Load-displacement
curve from the analytical model of the beam versus experimental data
from~\cite{bresler1963shear}.}%
\label{figure.loaddeflection}%
\end{figure}

We applied the proposed method to a simply supported
reinforced concrete beam. For the model problem we used specimen A-1 tested
in~\cite{bresler1963shear}, which is a simply supported beam subject to a
concentrated load at the center. A~detailed drawing is shown in
Figure~\ref{figure.beama1}. Due to the symmetry, only half of the beam is
considered in this study. The spatial discretization uses a two-dimensional
mesh with $60$ finite elements and $160$ degree of freedom under plane stress
conditions. Material properties of the beam are set as follows: the maximum
compressive stress of the concrete $\acute{f}_{c}=24.1\,$MPa, yield stress of
the bottom longitudinal reinforcement $f_{sy1}=555\,$MPa, yield stress of the
top longitudinal reinforcement $f_{sy2}=345\,$MPa and yield stress of the
stirrups $f_{sy3}=325\,$MPa. First, we implemented the deterministic
nonlinear finite element procedure as discussed in Section~\ref{sec:dfem} in
\textsc{Matlab}. A comparison of the load-\textcolor{black}{vertical displacement} curves from the
deterministic model with the experiment~\cite{bresler1963shear} can be seen in
Figure~\ref{figure.loaddeflection}. A good agreement between numerical and
experimental results can be seen throughout the entire load-\textcolor{black}{displacement} range.

In order to set up a stochastic model, we considered the initial concrete
modulus of elasticity, which is one of the parameters that have significant
effect on the response of a reinforced concrete member. There are many
parametrizations proposed in literature. In this study, we used a formula from
the FIB model code~\cite{code2010fib}, which is
\begin{equation}
E_{ci}=E_{c0}\alpha_{E}\left(  \frac{{f}_{cm}}{10}\right)  ^{\frac{1}{3}},
\label{equation.eci1}%
\end{equation}
where $E_{ci}$ is the initial modulus of elasticity in MPa at the concrete age
of $28$ days, $E_{c0}=21.5\,$GPa, $\alpha_{E}$ is a parameter that depend on
the types of aggregate, for example equal $0.7$ for Sandstone aggregate and
$0.9$ for Limestone aggregate, ${f}_{cm}$ is the actual compressive strength
of concrete at an age of $28$ days in MPa. In formula~(\ref{equation.eci1}),
we considered$~\alpha_{E}$ as a random input parameter with uniform
distribution and the mean value $0.8$. Two coefficients of variation (CoV) are
considered for the random input parameter that equal $5.77\%$ and $10\%$. That
is, the value of$~\alpha_{E}$ is uniformly distributed in the intervals
$(0.72-0.88)$ and $(0.66-0.94)$, respectively. The initial concrete modulus of
elasticity thus becomes stochastic, that is $E_{ci}\equiv E_{ci}(\xi)$.

Our goal in the setup is to find the stochastic expansion of the tangent
stiffness matrix~(\ref{equation.stiffdisc}) via~(\ref{equation.knj1}). This
calculation uses an expansion of the the tangent modulus of elasticity and is
repeated after every load increment. Based on the uniaxial stress-strain
relationship in compression, see eq.~(\ref{equation.conccomp}), the tangent
concrete modulus of elasticity is
\begin{equation}
E_{cT}=\frac{d\sigma_{c}}{d\epsilon_{c}}=E_{c1}\frac{(k-2\eta)-(k-2)\eta^{2}%
}{(1+(k-2)\eta)^{2}}. \label{equation.ect1}%
\end{equation}
Because $E_{cT}$ depends on $E_{ci}$, they are both stochastic, and in
particular
\begin{equation}
E_{cT}(\xi)=E_{c1}\frac{\left(  \frac{E_{ci}\left(  \xi\right)  }{E_{c1}%
}-2\frac{\epsilon_{c}}{\epsilon_{c1}}\right)  -\left(  \frac{E_{ci}(\xi
)}{E_{c1}}-2\right)  \left(  \frac{\epsilon_{c}}{\epsilon_{c1}}\right)  ^{2}%
}{\left(  1+\left(  \frac{E_{ci}(\xi)}{E_{c1}}-2\right)  \frac{\epsilon_{c}%
}{\epsilon_{c1}}\right)  ^{2}}. \label{equation.ect2}%
\end{equation}
Then, using also the values of$~\epsilon_{c}$ corresponding to the $n$th load
increment, we project~(\ref{equation.ect2}) on the gPC\ basis using numerical
quadrature to\ calculate the expansion of the tangent concrete modulus of
elasticity at the $n$th load increment
\begin{equation}
E_{cT}^{n}(\xi)=\sum_{i=1}^{n_{K}}E_{cTi}^{n}\psi_{i}(\xi).
\label{equation.enct}%
\end{equation}
Expansion~(\ref{equation.enct}) is then used to calculate
expansion~(\ref{equation.stiffdisc}). In general, while
by~\cite{Matthies-2005-GML} when considering the expansion of the
displacement~(\ref{equation.ui2}) using gPC\ polynomials of degree$~p$ it
would be possible to use gPC\ polynomials of degree up to$~2p$ in
expansion~(\ref{equation.enct}), in the numerical experiments we used
expansions of the same degree and we set $n_{K}=n_{\xi}$. In particular, we
did not observe with $n_{K}>n_{\xi}$ any improvement of results in our
numerical experiments.

In the numerical experiments, we show estimated probability density functions
(PDFs) of the$~E_{cT}$, and then we particularly focus on the vertical
displacement at the center of the beam$~u$, for which we also tabulate several
other statistical estimates. We note that$~E_{cT}$ is pertinent to one
element. We also compare the results of the stochastic Galerkin method (SG) to
those obtained by the stochastic collocation (SC) and Monte Carlo (MC) methods.
For the two spectral stochastic finite element methods, SC\ and SG, we used
gPC polynomials of degree four and eight and Smolyak sparse grid. The results
of the MC simulation are based on $10^{6}$ samples. As the validation
criterion, we used the root mean square error (RMSE) of the vertical
displacement at the center of the beam$~u$, defined as
\begin{equation}
\text{RMSE}=\sqrt{\frac{1}{n_{MC}}\sum_{i=1}^{n_{MC}}\left(  u_{\star}\left(
\xi^{(i)}\right)  -u_{MC}\left(  \xi^{(i)}\right)  \right)  ^{2}},
\label{eq:RMSE}%
\end{equation}
where $\star$ indicates either SC or SG\ method. We also report in tables
estimates of the mean~$\mu$, standard deviation~$\sigma$ and probability
\textcolor{black}{${Pr}(u-k u_{m}\geq0)$, where $u_{m}$ is the vertical displacement 
at the center of the beam corresponding to the deterministic problem (the parameters are set
to the mean values), and values of~$k$ are specified in the tables.}

\begin{figure}[b]
\centering
\begin{tabular}
[c]{cc}%
\subfigure[$E_{cT}$ at load\ increment~10 \label{figure:ect57710}] {\includegraphics[width=8cm]{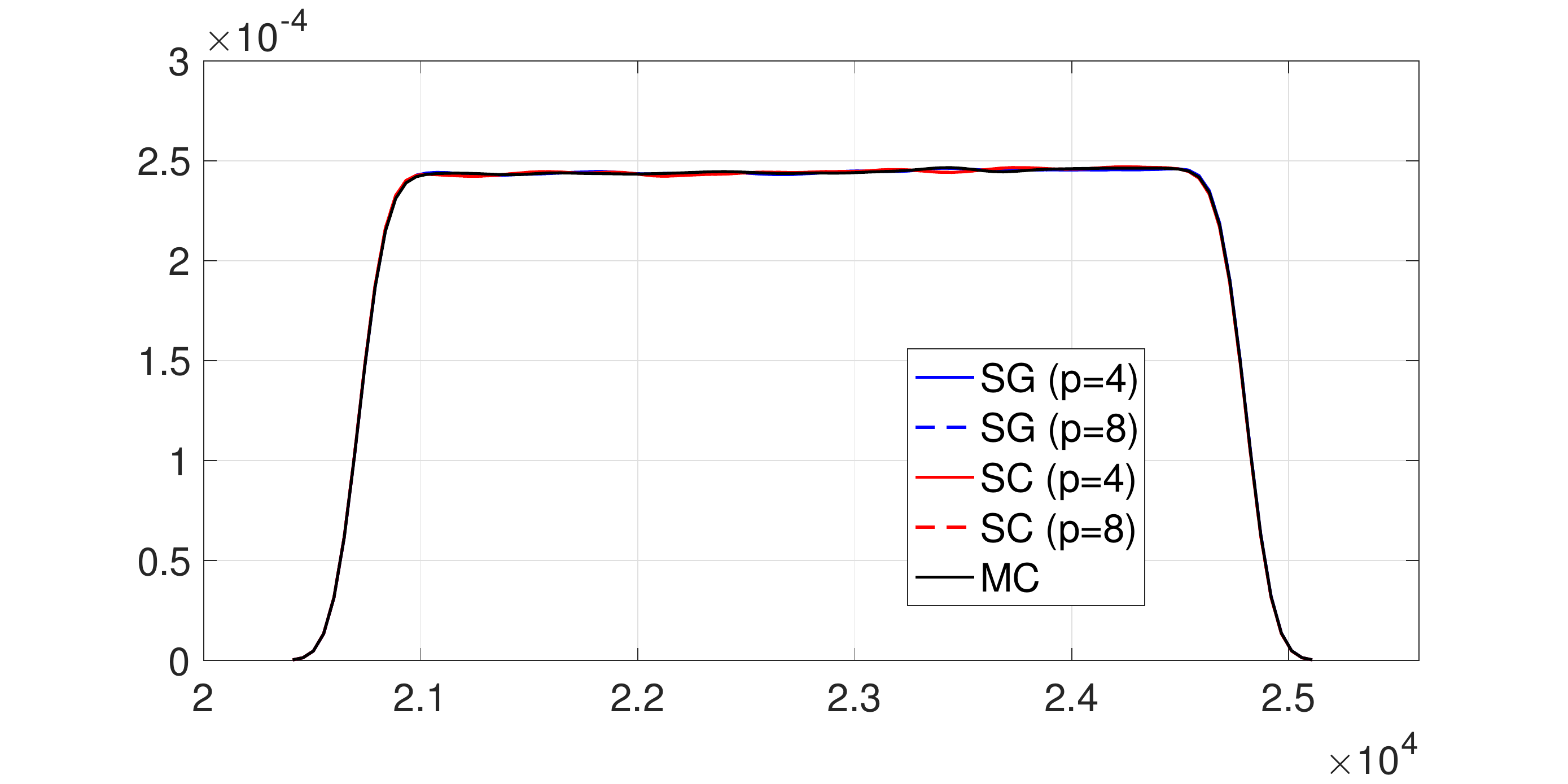}} &
\subfigure[displacement at load\ increment~10 \label{figure:u57710}] {\includegraphics[width=8cm]{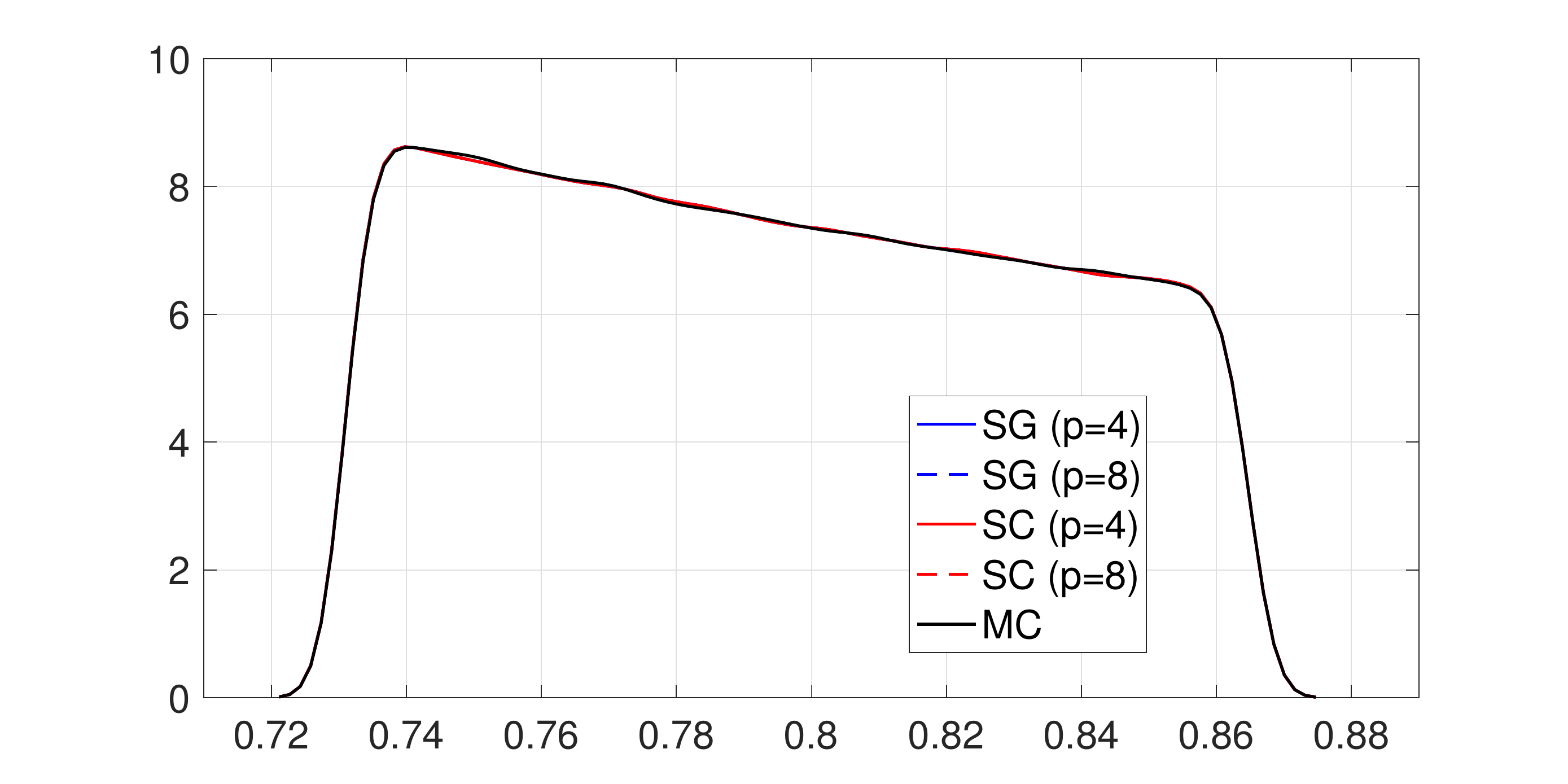}}\\
\subfigure[$E_{cT}$ at load\ increment\ 20 \label{figure:ect57720}] {\includegraphics[width=8cm]{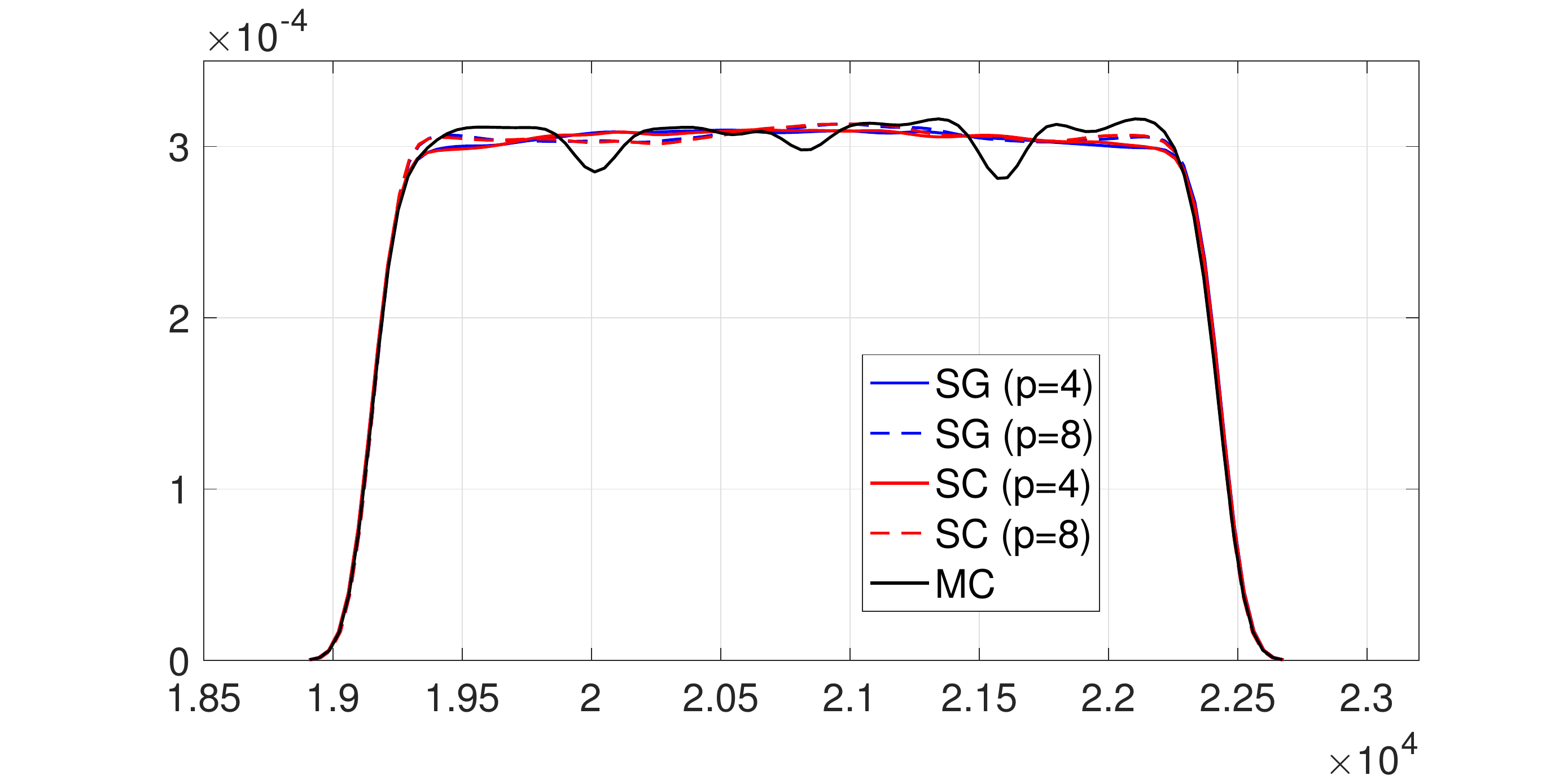}} &
\subfigure[displacement at load\ increment\ 20 \label{figure:u57720}] {\includegraphics[width=8cm]{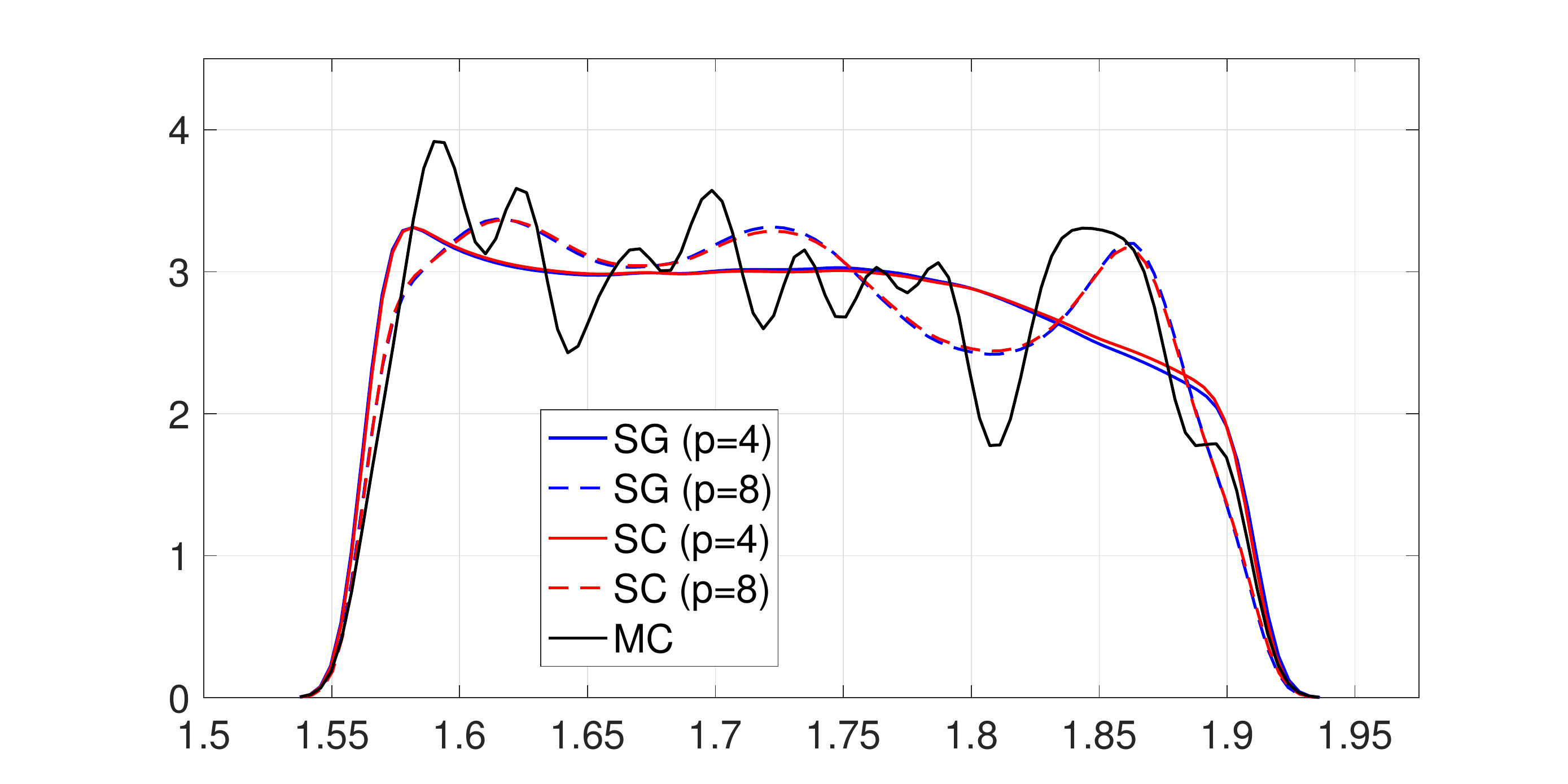}}\\
\subfigure[$E_{cT}$ at load\ increment\ 87 \label{figure:ect57787}] {\includegraphics[width=8cm]{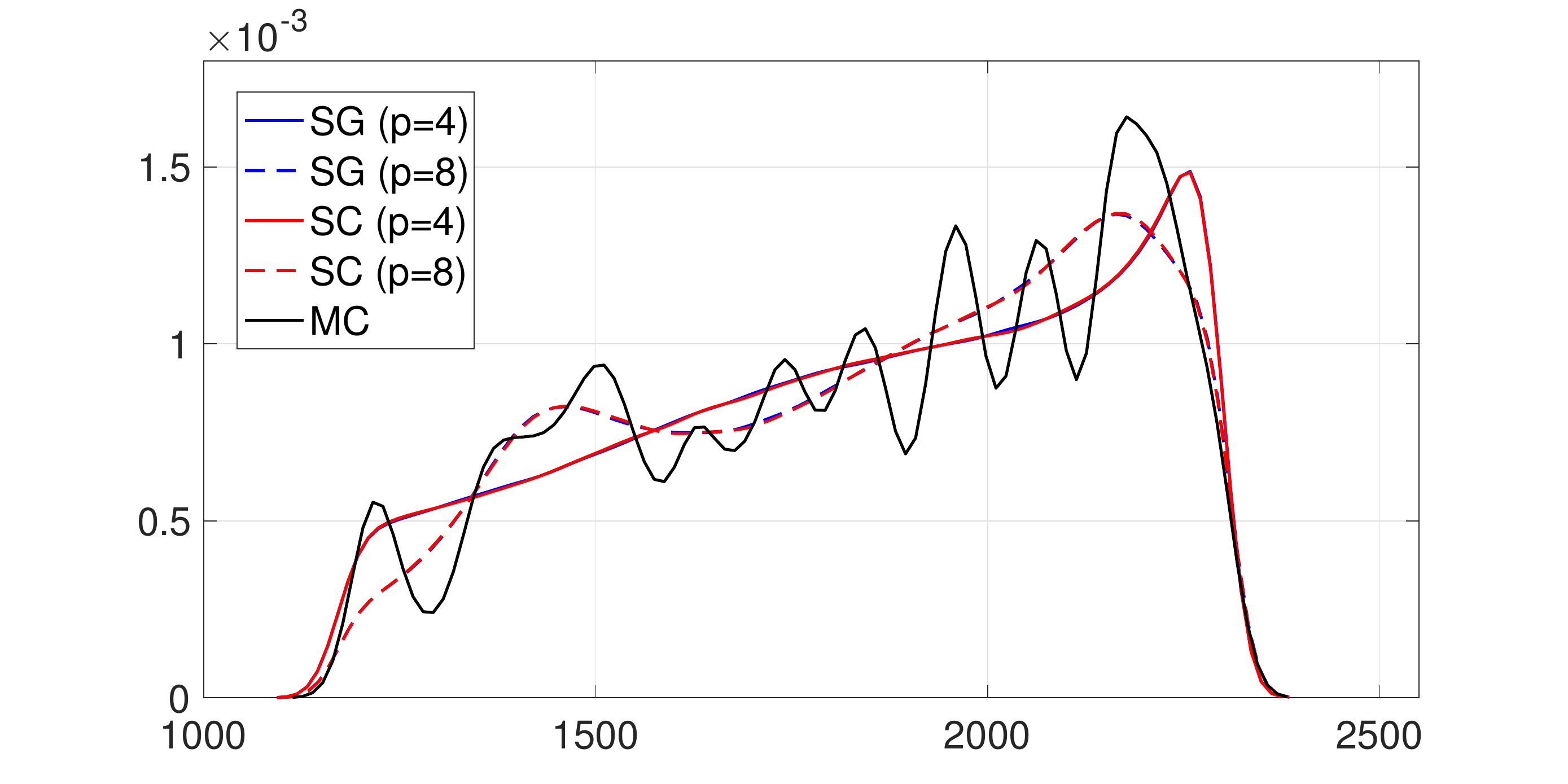}} &
\subfigure[displacement at load\ increment\ 87 \label{figure:u57787}] {\includegraphics[width=8cm]{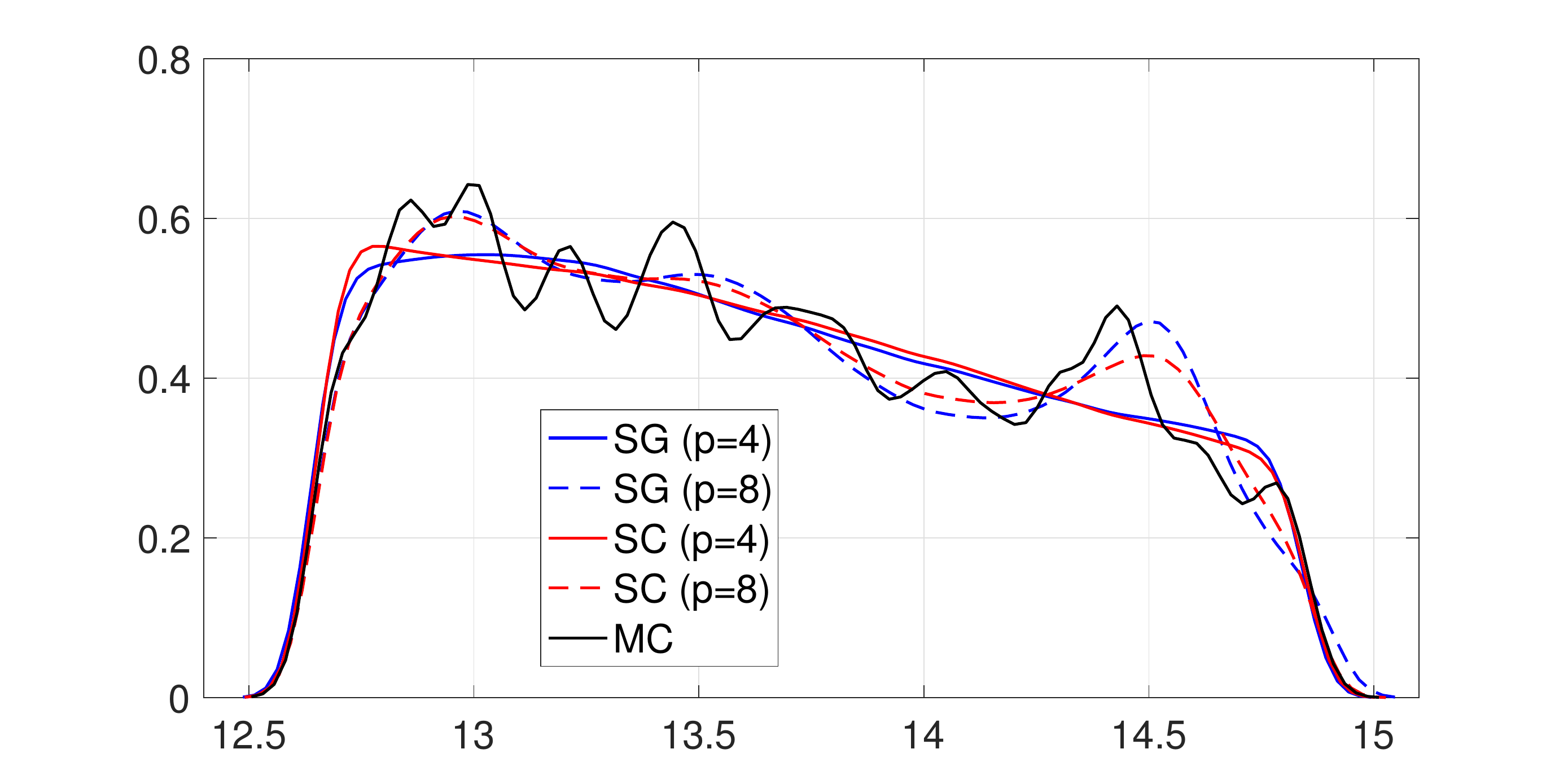}}\\
&
\end{tabular}
\caption{Estimated PDFs of $E_{cT}$ and the displacement at center of the beam
during the loading, $CoV=5.77\%$.}%
\label{figure:CoV577}%
\end{figure}

Figure~\ref{figure:CoV577} shows the estimated PDFs of the$~E_{cT}$ and the
vertical displacement at the center of the beam for the case with $CoV=5.77\%$
at load increments $10$, $20$ and $87$, which is the final load increment in
this case. A crack starts to propagate between the load increments $10$ and
$20$ (specifically at the increment$~14$). Correspondingly, it can be observed
that at the load increment$~$10 all PDF estimates are smooth and coincide, but
at the increment$~20$ the PDF\ obtained by MC\ is oscillatory and the SC\ and
SG\ methods provide a smooth interpolant. We also observe that the
gPC\ polynomial with degree $p=8$ is relatively more oscillatory than the gPC
polynomial with degree $p=4$. From the plots we anticipate that a polynomial
of a very high degree, unfeasible for practical computations, would be needed
in order to match the results of the MC\ simulation. Nevertheless, it turns
out that even the low-degree gPC approximation provides an accurate
quantitative insight into the response of the beam to loading. Such insight is
provided by Table~\ref{tab:1}, which compares estimated values of the
mean~$\mu$ and standard deviation~$\sigma$ for the displacement at center of
the beam with $CoV=5.77\%$ using Monte Carlo (MC), stochastic Galerkin (SG)
and stochastic collocation (SC) methods with gPC degrees $p=4$ and $8$, and an
assessment of the methods using both the RMSE value~(\ref{eq:RMSE}).
\textcolor{black}{Table~\ref{tab:1a} then provides estimates of
probability $\operatorname{Pr}(u-\widetilde{u}\geq0)$, where $u$ is the displacement
at the center of the beam from the stochastic problem and different setting of $\widetilde{u}$, in which 
$u_{m}$ is the displacement at the center of the beam from the deterministic problem with the
mean value parameters.} 
We see a good agreement of all methods in \textcolor{black}{all} indicators
$\mu$, $\sigma$ and $\Pr$. In fact, from the results we cannot discern any
improvement by increasing the gPC\ degree. 
 
\begin{figure}[b]
\centering
\begin{tabular}
[c]{cc}%
\subfigure[$E_{cT}$ at load\ increment\ 10 \label{figure:ect1010}] {\includegraphics[width=8cm,height=5.2cm]{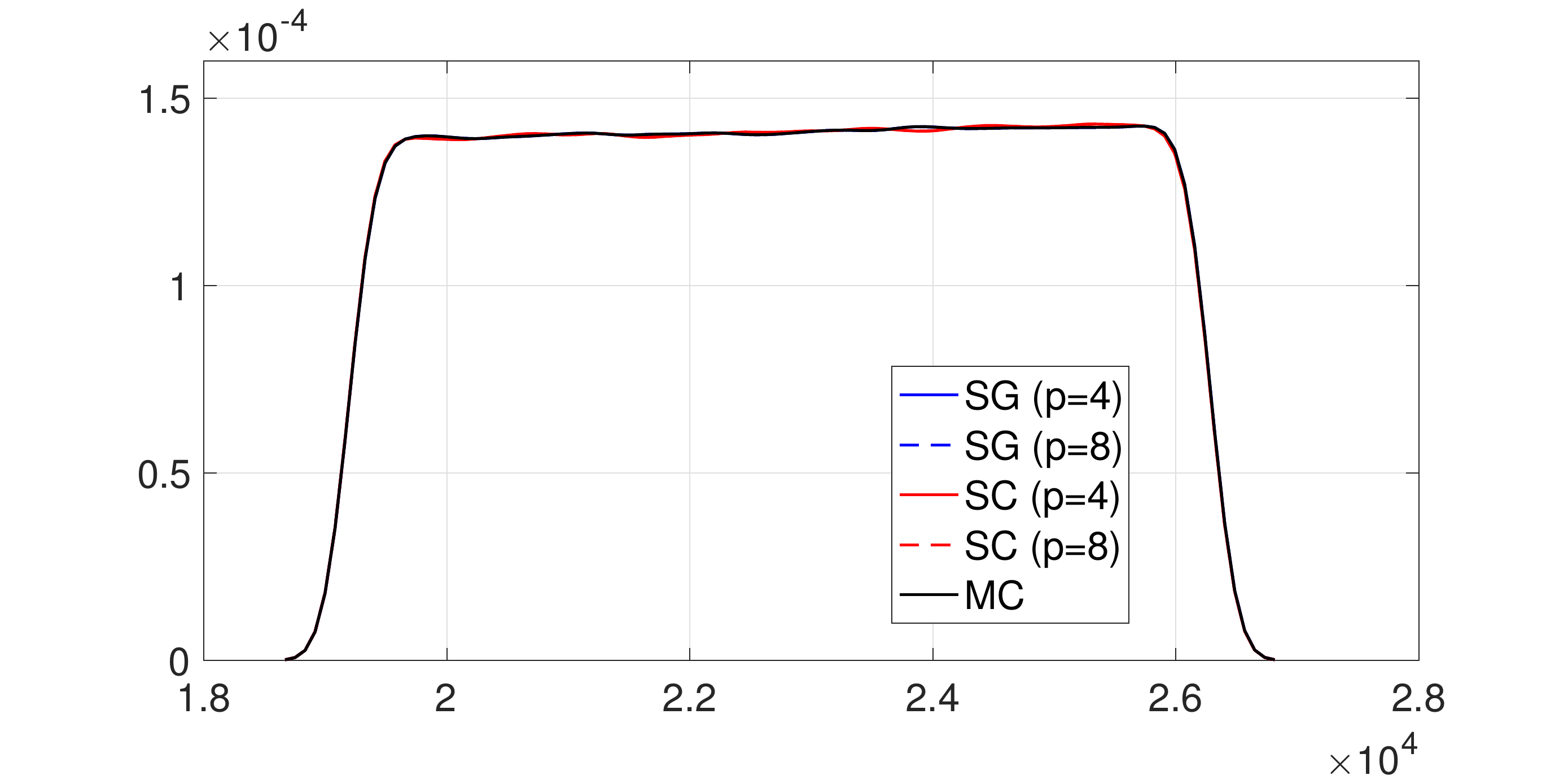}} &
\subfigure[displacement at load\ increment\ 10 \label{figure:u1010}] {\includegraphics[width=8cm,height=5.2cm]{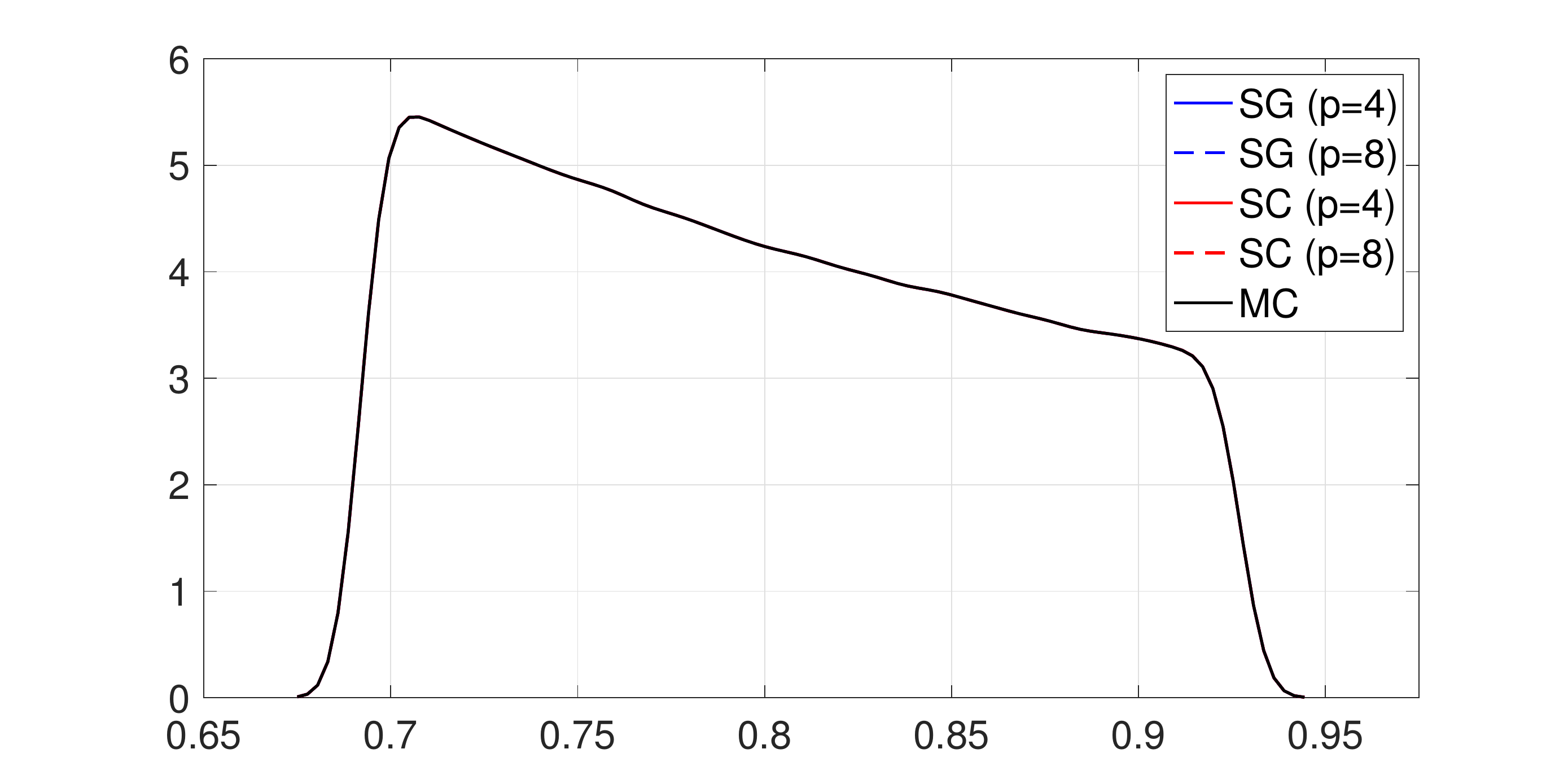}}\\
\subfigure[$E_{cT}$ at load\ increment\ 20 \label{figure:ect1020}] {\includegraphics[width=8cm,height=5.2cm]{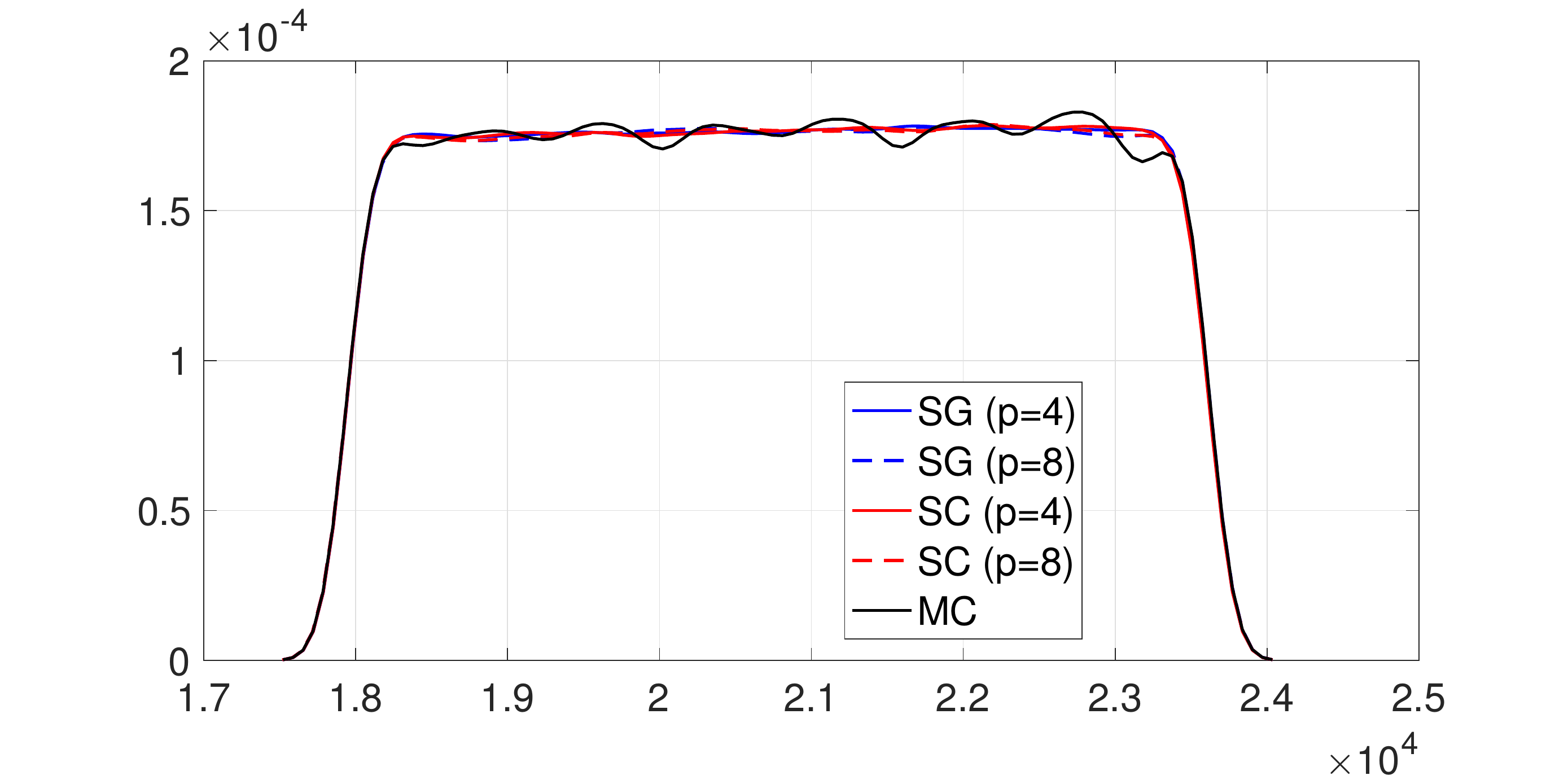}} &
\subfigure[displacement at load\ increment\ 20 \label{figure:u1020}] {\includegraphics[width=8cm,height=5.2cm]{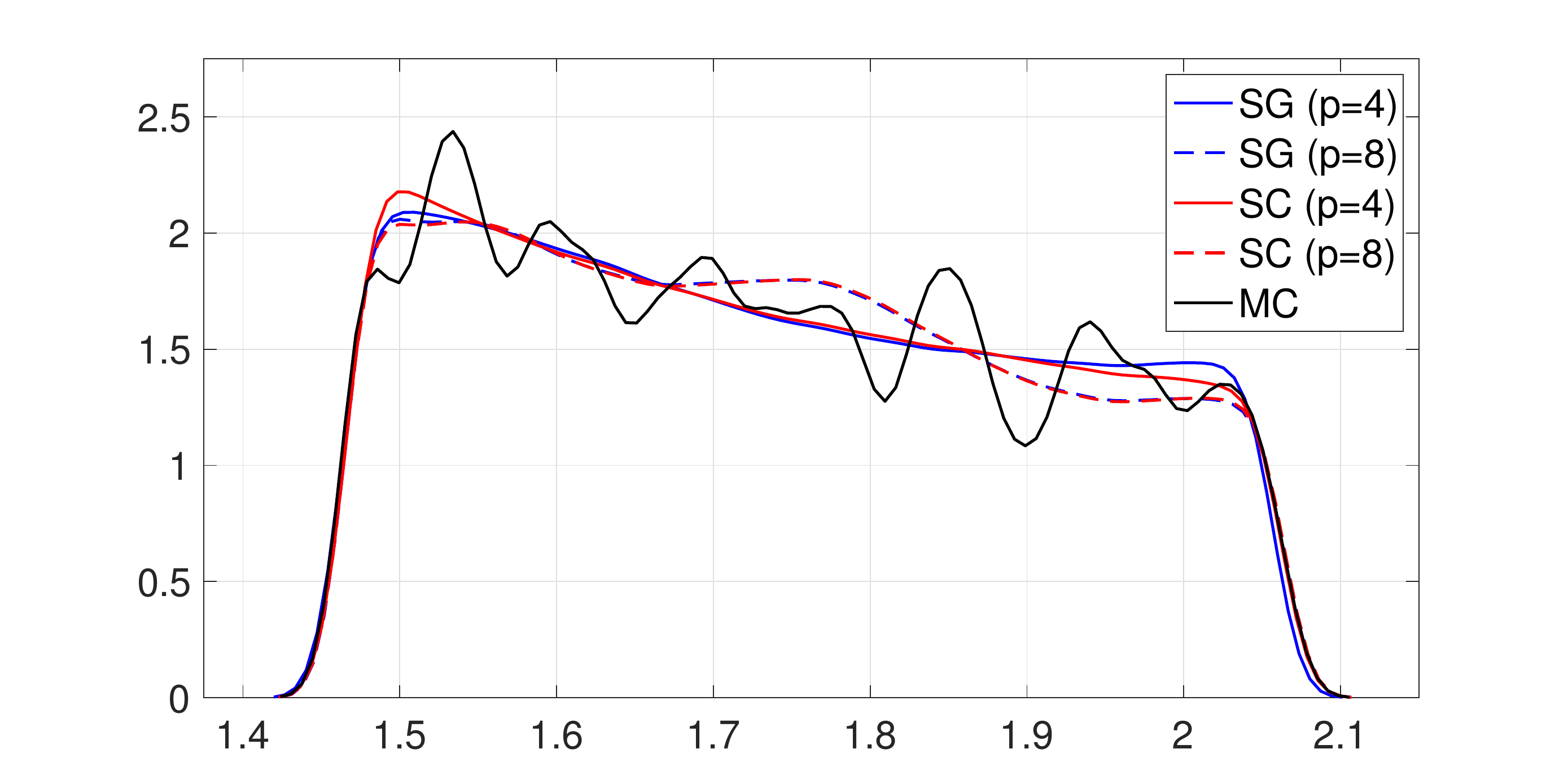}}\\
\subfigure[$E_{cT}$ at load\ increment\ 86 \label{figure:ect1086}] {\includegraphics[width=8cm,height=5.2cm]{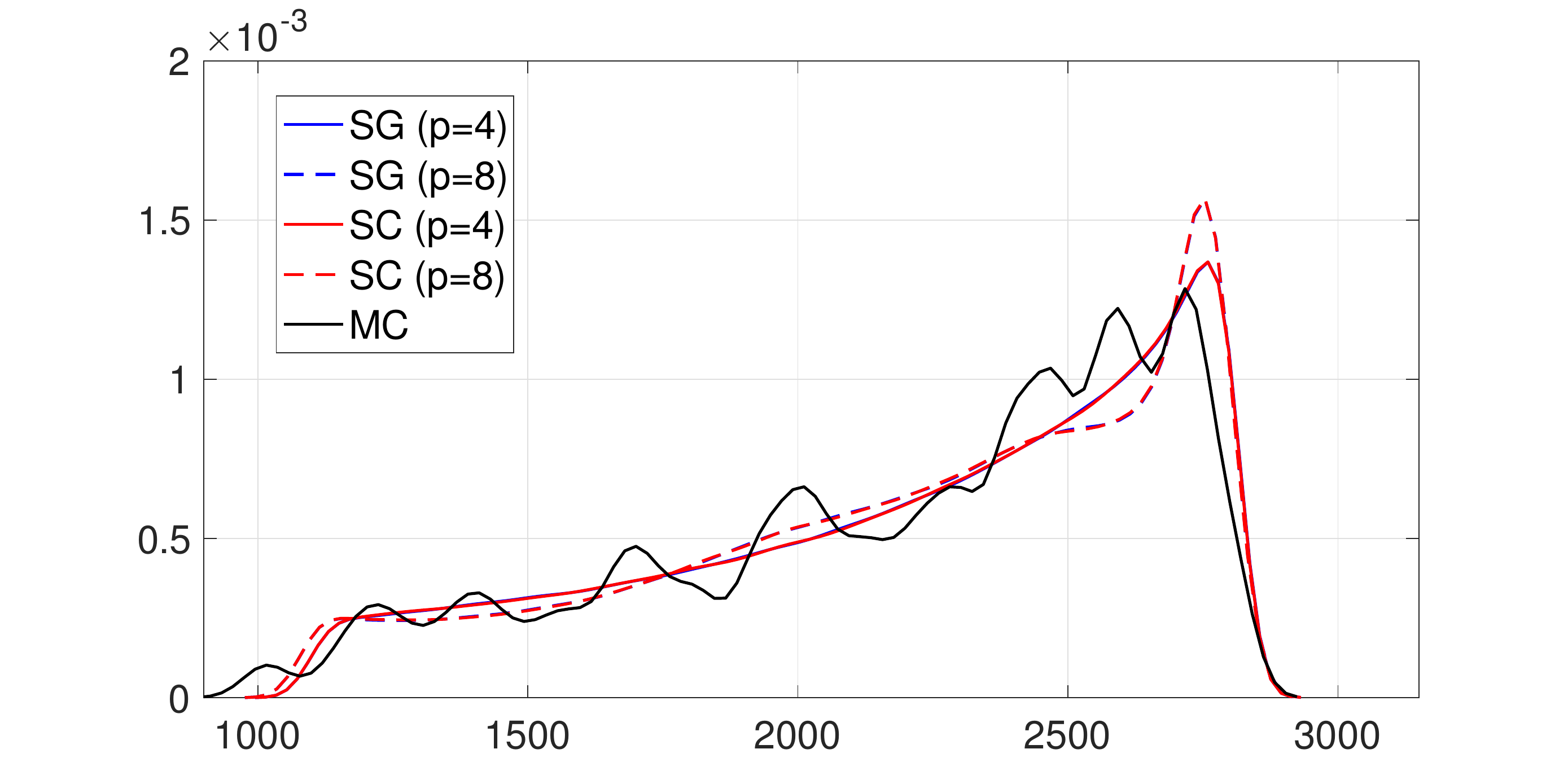}} &
\subfigure[displacement at load\ increment\ 86 \label{figure:u1086}] {\includegraphics[width=8cm,height=5.2cm]{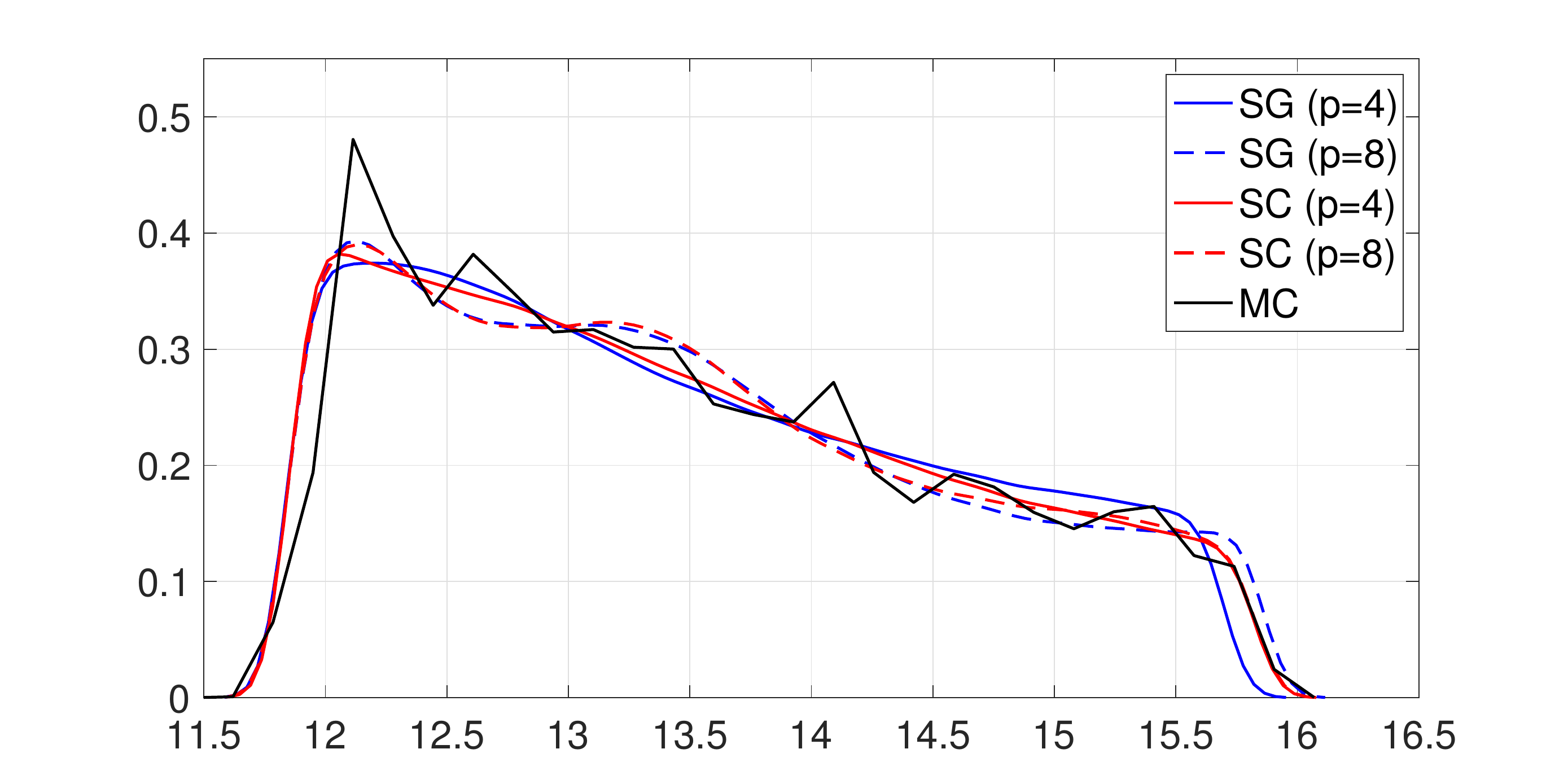}}\\
&
\end{tabular}
\caption{Estimated PDFs of $E_{cT}$ and the displacement at center of the beam
during the loading, $CoV=10\%$.}%
\label{figure:CoV10}%
\end{figure}

Figure~\ref{figure:CoV10}, Table~\ref{tab:2} \textcolor{black}{and Table~\ref{tab:2a}} then show results
corresponding to $CoV=10\%$. Comparing to the case with $CoV=5.77\%$, we see
that all PDFs have larger support, but the qualitative behavior of the
solutions is similar to the previous case. In particular, at the load
increment$~10$, the results of all methods match, and after the crack starts
to propagate at higher load increments, the gPC\ methods provide a smooth
interpolation to the MC\ solution. Nevertheless, in all cases the statistical
indicators $\mu$, $\sigma$ and $\Pr$ reported in Table~\ref{tab:2} \textcolor{black}{and Table~\ref{tab:2a}} are \textcolor{black}{again} in a
good agreement with Monte Carlo simulation.

\begin{table}[b]
\caption{Displacement at center of the beam with $CoV=5.77\%$: estimated
values of the mean~$\mu$ and standard deviation~$\sigma$ by Monte Carlo (MC),
stochastic Galerkin (SG) and stochastic collocation (SC) methods with gPC
degrees $p=4$ and $8$, comparison of the methods using root mean square error
(RMSE), see eq.~(\ref{eq:RMSE}).
}%
\label{tab:1}
\begin{center}%
\begin{tabular}
[c]{|c|c|c|c|}\hline
method & $\mu$ & $\sigma$ & $\operatorname{RMSE}$ \\\hline
\multicolumn{4}{|c|}{load increment 10}\\\hline
MC & $7.9414\times10^{-1}$ & $3.8566\times10^{-2}$ & - \\\hline
SG, $p=4$ & $7.9415\times10^{-1}$ & $3.8567\times10^{-2}$ & $4.7201\times
10^{-2}$ \\
SG, $p=8$ & $7.9415\times10^{-1}$ & $3.8567\times10^{-2}$ & $4.7201\times
10^{-2}$ \\\hline
SC, $p=4$ & $7.9414\times10^{-1}$ & $3.8566\times10^{-2}$ & $4.7201\times
10^{-2}$ \\
SC, $p=8$ & $7.9415\times10^{-1}$ & $3.8567\times10^{-2}$ & $4.7201\times
10^{-2}$ \\\hline
\multicolumn{4}{|c|}{load increment 20}\\\hline
MC & $1.7267$ & $9.8225\times10^{-2}$ & - \\\hline
SG, $p=4$ & $1.7261$ & $9.8658\times10^{-2}$ & $1.2050\times10^{-1}$ \\
SG, $p=8$ & $1.7260$ & $9.7058\times10^{-2}$ & $1.1952\times10^{-1}$ \\\hline
SC, $p=4$ & $1.7263$ & $9.8590\times10^{-2}$ & $1.2046\times10^{-1}$ \\
SC, $p=8$ & $1.7260$ & $9.7035\times10^{-2}$ & $1.1951\times10^{-1}$ \\\hline
\multicolumn{4}{|c|}{load increment 87}\\\hline
MC & $13.6311$ & $6.1582\times10^{-1}$ & - \\\hline
SG, $p=4$ & $13.6280$ & $6.1808\times10^{-1}$ & $7.5516\times10^{-1}$ \\
SG, $p=8$ & $13.6401$ & $6.2000\times10^{-1}$ & $7.5639\times10^{-1}$ \\\hline
SC, $p=4$ & $13.6275$ & $6.1693\times10^{-1}$ & $7.5445\times10^{-1}$ \\
SC, $p=8$ & $13.6309$ & $6.1700\times10^{-1}$ & $7.5449\times10^{-1}$ \\\hline
\end{tabular}
\end{center}
\end{table}

\begin{table}[b]
\caption{Estimated probability $\operatorname{Pr}%
(u-\widetilde{u}\geq0)$ in \%, where $u$ is the displacement at the center of the beam from
the stochastic problem with $CoV=5.77\%$ and different setting of $\widetilde{u}$, in which $u_{m}$ is the displacement at the center of the
beam from the deterministic problem with the mean value parameters.}%
\label{tab:1a}
\begin{center}%
\begin{tabular}
[c]{|c|c|c|c|c|c|c|}\hline
method & $\widetilde{u}=u_{m}$ & $\widetilde{u}=1.02u_{m}$ & $\widetilde{u}=1.04u_{m}$ & $\widetilde{u}=1.06u_{m}$ &$\widetilde{u}=1.08$ &$\widetilde{u}=1.10u_{m}$\\\hline
\multicolumn{7}{|c|}{load increment 10}\\\hline
MC & $48.56\,\%$ & $36.93 \,\%$ & $25.75 \,\%$ & $14.97 \,\%$ & $4.56 \,\%$ & $0 \,\%$\\\hline
SG, $p=4$ & $48.57\,\%$ & $36.95 \,\%$ & $25.74 \,\%$ & $14.95 \,\%$ & $4.55 \,\%$ & $0 \,\%$\\
SG, $p=8$ & $48.57\,\%$ & $36.95 \,\%$ & $25.74 \,\%$ & $14.95 \,\%$ & $4.55 \,\%$ & $0 \,\%$\\\hline
SC, $p=4$ & $48.57\,\%$ & $36.95 \,\%$ & $25.74 \,\%$ & $14.95 \,\%$ & $4.55 \,\%$ & $0 \,\%$\\
SC, $p=8$ & $48.57\,\%$ & $36.95 \,\%$ & $25.74 \,\%$ & $14.95 \,\%$ & $4.55 \,\%$ & $0 \,\%$\\\hline
\multicolumn{7}{|c|}{load increment 20}\\\hline
MC & $49.22\,\%$ & $38.79 \,\%$ & $28.11 \,\%$ & $20.98 \,\%$ & $9.57 \,\%$ & $1.87 \,\%$\\\hline
SG, $p=4$ & $49.17\,\%$ & $38.73 \,\%$ & $28.49 \,\%$ & $18.83 \,\%$ & $10.14 \,\%$ & $2.44 \,\%$\\
SG, $p=8$ & $48.24\,\%$ & $37.34 \,\%$ & $28.34 \,\%$ & $19.94 \,\%$ & $9.70 \,\%$ & $1.44 \,\%$\\\hline
SC, $p=4$ & $49.27\,\%$ & $38.84 \,\%$ & $28.64 \,\%$ & $18.97 \,\%$ & $10.15 \,\%$ & $2.25 \,\%$\\
SC, $p=8$ & $48.30\,\%$ & $37.43 \,\%$ & $28.35 \,\%$ & $19.88 \,\%$ & $9.67 \,\%$ & $1.51 \,\%$\\\hline
\multicolumn{7}{|c|}{load increment 87}\\\hline
MC & $47.09\,\%$ & $34.33 \,\%$ & $24.16 \,\%$ & $12.54 \,\%$ & $3.81 \,\%$ & $0 \,\%$\\\hline
SG, $p=4$ & $46.98\,\%$ & $34.50 \,\%$ & $23.24 \,\%$ & $13.12 \,\%$ & $3.84 \,\%$ & $0 \,\%$\\
SG, $p=8$ & $46.87\,\%$ & $34.68 \,\%$ & $24.90 \,\%$ & $14.35 \,\%$ & $3.34 \,\%$ & $0 \,\%$\\\hline
SC, $p=4$ & $47.21\,\%$ & $34.50 \,\%$ & $23.04 \,\%$ & $12.89 \,\%$ & $3.84 \,\%$ & $0 \,\%$\\
SC, $p=8$ & $46.83\,\%$ & $34.71 \,\%$ & $24.39 \,\%$ & $13.77 \,\%$ & $3.25 \,\%$ & $0 \,\%$\\\hline
\end{tabular}
\end{center}
\end{table}

\begin{table}[b]
\caption{Displacement at center of the beam with $CoV=10\%$: the headers are
the same as in Table~\ref{tab:1}.}%
\label{tab:2}
\begin{center}%
\begin{tabular}
[c]{|c|c|c|c|}\hline
& $\mu$ & $\sigma$ & $\operatorname{RMSE}$ \\\hline
\multicolumn{4}{|c|}{load increment 10}\\\hline
MC & $7.9804\times10^{-1}$ & $6.7548\times10^{-2}$ & - \\\hline
SG, $p=4$ & $7.9804\times10^{-1}$ & $6.7548\times10^{-2}$ & $2.5263\times
10^{-6}$ \\
SG, $p=8$ & $7.9804\times10^{-1}$ & $6.7548\times10^{-2}$ & $0.0332\times
10^{-6}$ \\\hline
SC, $p=4$ & $7.9804\times10^{-1}$ & $6.7548\times10^{-2}$ & $2.5284\times
10^{-6}$ \\
SC, $p=8$ & $7.9804\times10^{-1}$ & $6.7548\times10^{-2}$ & $0.0331\times
10^{-6}$ \\\hline
\multicolumn{4}{|c|}{load increment 20}\\\hline
MC & $1.7365$ & $1.7224\times10^{-1}$ & - \\\hline
SG, $p=4$ & $1.7367$ & $1.7247\times10^{-1}$ & $3.3543\times10^{-3}$ \\
SG, $p=8$ & $1.7363$ & $1.7242\times10^{-1}$ & $3.4802\times10^{-3}$ \\\hline
SC, $p=4$ & $1.7365$ & $1.7255\times10^{-1}$ & $3.0136\times10^{-3}$ \\
SC, $p=8$ & $1.7365$ & $1.7247\times10^{-1}$ & $3.3217\times10^{-3}$ \\\hline
\multicolumn{4}{|c|}{load increment 86}\\\hline
MC & $13.4660$ & $1.0815\times10^{0}$ & - \\\hline
SG, $p=4$ & $13.4496$ & $1.0789\times10^{0}$ & $5.6684\times10^{-2}$ \\
SG, $p=8$ & $13.4615$ & $1.1007\times10^{0}$ & $5.2660\times10^{-2}$ \\\hline
SC, $p=4$ & $13.4568$ & $1.0911\times10^{0}$ & $4.9047\times10^{-2}$ \\
SC, $p=8$ & $13.4565$ & $1.0907\times10^{0}$ & $5.0035\times10^{-2}$ \\\hline
\end{tabular}
\end{center}
\end{table}

\begin{table}[b]
\caption{Estimated probability $\operatorname{Pr}%
(u-\widetilde{u}\geq0)$ in \%, where $u$ is the displacement at the center of the beam from
the stochastic problem with $CoV=10\%$ and different setting of $\widetilde{u}$, in which $u_{m}$ is the displacement at the center of the
beam from the deterministic problem with the mean value parameters.}%
\label{tab:2a}
\begin{center}%
\begin{tabular}
[c]{|c|c|c|c|c|c|c|}\hline
method & $\widetilde{u}=u_{m}$ & $\widetilde{u}=1.04u_{m}$ & $\widetilde{u}=1.08u_{m}$ & $\widetilde{u}=1.12u_{m}$ &$\widetilde{u}=1.16u_{m}$ &$\widetilde{u}=1.20u_{m}$\\\hline
\multicolumn{7}{|c|}{load increment 10}\\\hline
MC & $47.51\,\%$ & $34.40 \,\%$ & $22.26 \,\%$ & $11.00 \,\%$ & $0.52 \,\%$ & $0 \,\%$\\\hline
SG, $p=4$ & $47.51\,\%$ & $34.40 \,\%$ & $22.26 \,\%$ & $11.00 \,\%$ & $0.52 \,\%$ & $0 \,\%$\\
SG, $p=8$ & $47.51\,\%$ & $34.40 \,\%$ & $22.26 \,\%$ & $11.00 \,\%$ & $0.52 \,\%$ & $0 \,\%$\\\hline
SC, $p=4$ & $47.51\,\%$ & $34.40 \,\%$ & $22.26 \,\%$ & $11.00 \,\%$ & $0.52 \,\%$ & $0 \,\%$\\
SC, $p=8$ & $47.51\,\%$ & $34.40 \,\%$ & $22.26 \,\%$ & $11.00 \,\%$ & $0.52 \,\%$ & $0 \,\%$\\\hline
\multicolumn{7}{|c|}{load increment 20}\\\hline
MC & $47.47\,\%$ & $36.35 \,\%$ & $24.72 \,\%$ & $15.96 \,\%$ & $6.52 \,\%$ & $0 \,\%$\\\hline
SG, $p=4$ & $47.68\,\%$ & $36.67 \,\%$ & $26.21 \,\%$ & $16.12 \,\%$ & $6.15 \,\%$ & $0 \,\%$\\
SG, $p=8$ & $47.96\,\%$ & $35.72 \,\%$ & $25.45 \,\%$ & $16.17 \,\%$ & $6.29 \,\%$ & $0 \,\%$\\\hline
SC, $p=4$ & $47.64\,\%$ & $36.53 \,\%$ & $25.98 \,\%$ & $15.97 \,\%$ & $6.37 \,\%$ & $0 \,\%$\\
SC, $p=8$ & $47.99\,\%$ & $35.91 \,\%$ & $25.47 \,\%$ & $16.08 \,\%$ & $6.46 \,\%$ & $0 \,\%$\\\hline
\multicolumn{7}{|c|}{load increment 86}\\\hline
MC & $44.69\,\%$ & $31.16 \,\%$ & $20.26 \,\%$ & $10.91 \,\%$ & $2.70 \,\%$ & $0 \,\%$\\\hline
SG, $p=4$ & $44.91\,\%$ & $31.53 \,\%$ & $20.09 \,\%$ & $10.13 \,\%$ & $1.20 \,\%$ & $0 \,\%$\\
SG, $p=8$ & $44.79\,\%$ & $30.48 \,\%$ & $19.82 \,\%$ & $11.26 \,\%$ & $3.47 \,\%$ & $0 \,\%$\\\hline
SC, $p=4$ & $44.94\,\%$ & $31.23 \,\%$ & $19.92 \,\%$ & $10.58 \,\%$ & $2.66 \,\%$ & $0 \,\%$\\
SC, $p=8$ & $44.70\,\%$ & $30.49 \,\%$ & $19.89 \,\%$ & $10.86 \,\%$ & $2.73 \,\%$ & $0 \,\%$\\\hline
\end{tabular}
\end{center}
\end{table}

{\color{black}{
\subsection{Example~2: Reinforced concrete shear wall}

\label{sec:sfem:example2} We also consider a reinforced concrete shear wall as specimen SH-L tested in~\cite{mickleborough1999prediction}. A~detailed drawing is shown in Figure~\ref{figure:shearwall}. The spatial discretization uses a two-dimensional mesh with $72$ finite elements and $182$ degree of freedom under plane stress conditions. Material properties of the shear wall are set as follows: the maximum compressive stress of the concrete $\acute{f}_{c}=44.7\,$MPa, yield stress of the horizontal and vertical reinforcement $f_{sy1}=460\,$MPa. The wall is first subjected to a vertical load of $357$KN, then the lateral load is applied incrementally. First, we implemented the deterministic nonlinear finite element procedure as discussed in Section~\ref{sec:sfem:galerkin} in \textsc{Matlab}. A comparison of the load-\textcolor{black}{horizontal displacement} curves from the deterministic model with the experiment~\cite{mickleborough1999prediction} can be seen in
Figure~\ref{figure:loaddeflection2}. A good agreement between numerical and experimental results can be seen throughout the entire load-\textcolor{black}{displacement} range.

\begin{figure}[b]
\centering
\includegraphics[width=12cm,height=7cm]{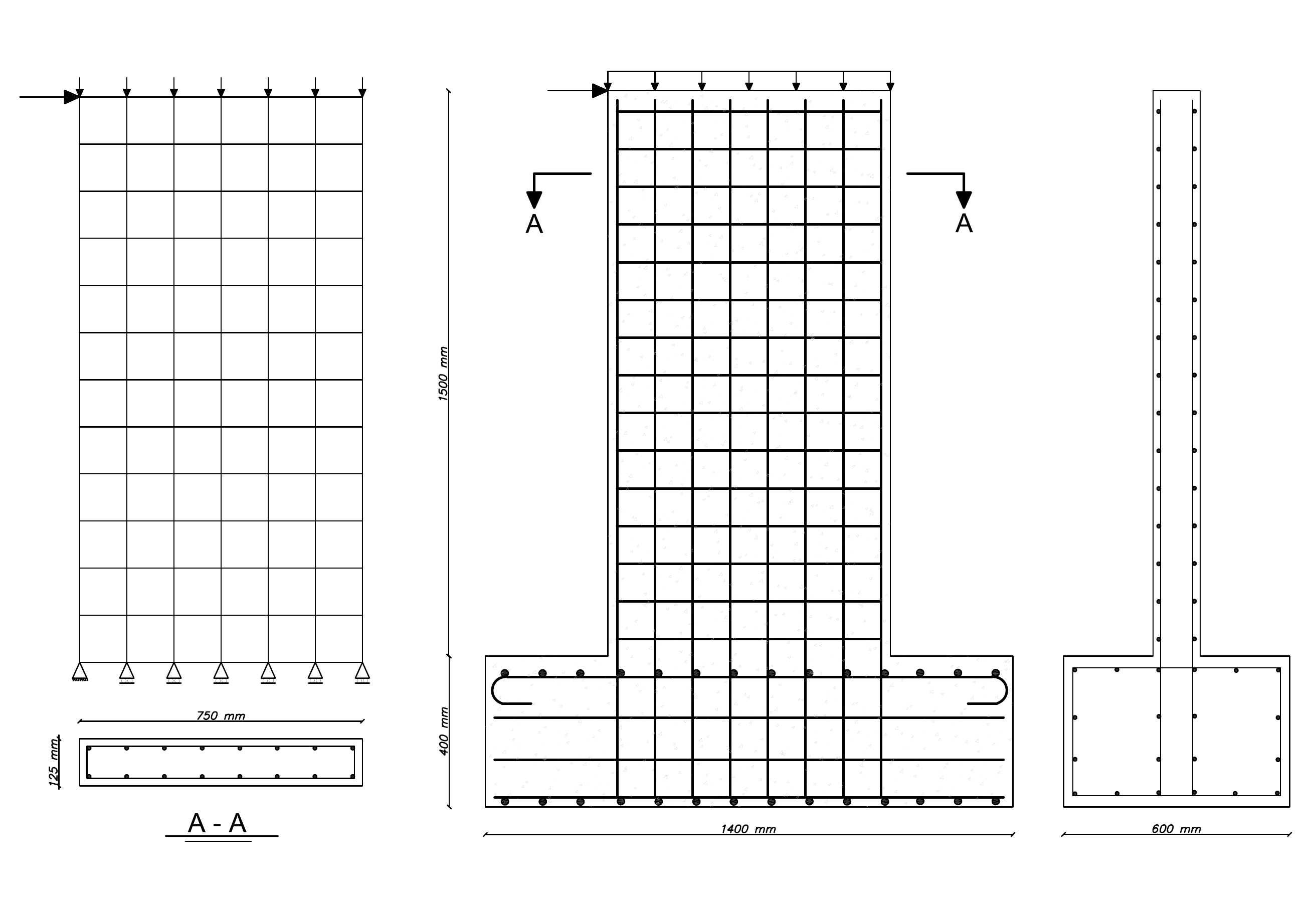}
\caption{Technical drawing of the shear wall. }%
\label{figure:shearwall}%
\end{figure}

\begin{figure}[b]
\centering
\includegraphics[width=12cm]{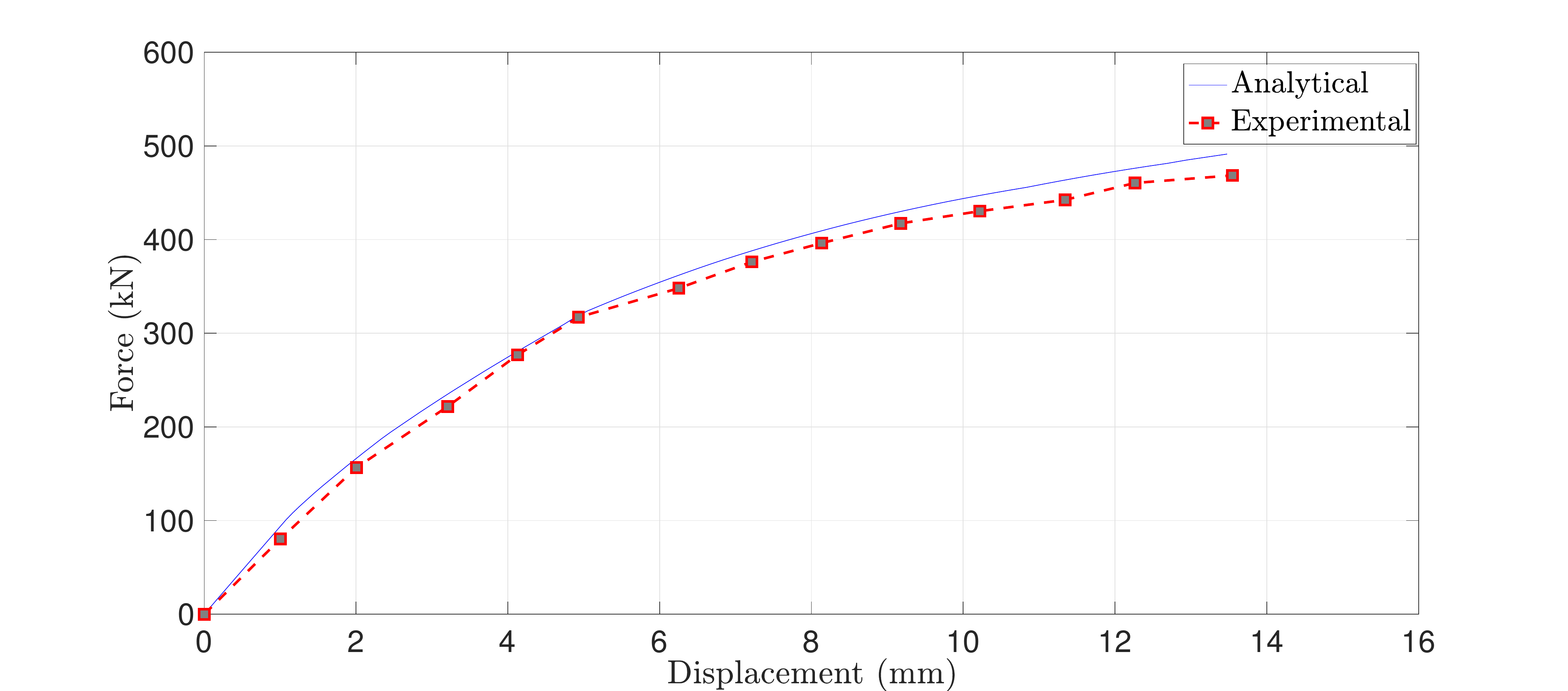} \caption{Load-displacement
curve from the analytical model of the shear wall versus experimental data
from~\cite{mickleborough1999prediction}.}%
\label{figure:loaddeflection2}%
\end{figure}

For the stochastic model, we considered the maximum compressive stress of the concrete $\acute{f}_{c}$, as a random input parameter with uniform distribution. To consider the spatial variability of the random parameter, we assume that mean of the maximum compressive stress of the concrete changes linearly from $(44.7+1.5)$MPa in the lowest mesh row to $(44.7-1.5)$MPa in the uppermost mesh row. The coefficients of variation (CoV) is considered equal to~$5\%$.

In the numerical experiments we show estimated probability density function (PDF) of the displacement~$u$ at the location of the applied lateral load,
and we also tabulate several other statistical estimates including RMSE from~(\ref{eq:RMSE}). We compare the results of the stochastic Galerkin method (SG) to those obtained by stochastic collocation (SC) and Monte Carlo (MC) methods.
For the two spectral stochastic finite element methods, SC\ and SG, we used gPC polynomials of degree~4 and~8 and Smolyak sparse grid. The results of the MC simulation are based on $10^{6}$ samples.

Figure~\ref{figure:CoV5} shows the estimated PDFs of the displacement at the location of the lateral load at top of the shear wall at load increments $5$, $10$, $20$ and $50$, which is the final load increment. 
A crack starts to propagate between the load increments~$10$ and~$20$ since that is when the PDF obtained by MC becomes oscillatory. 
It can be also seen, similarly as for the beam, that both SC and SG methods provide a smooth interpolant with degree $p=8$ relatively more oscillatory than with degree $p=4$.
Table~\ref{tab:3} then compares estimated values of the mean~$\mu$ and standard deviation~$\sigma$ for the displacement at the location of the lateral load at top of the shear wall using MC, SG and SC methods and an assessment of the methods using both the RMSE and probability estimates. 
It can be seen again that all methods are in quite a good agreement with Monte Carlo simulation, and this holds in particular for the probability estimates. 

\begin{figure}[b]
\centering
\begin{tabular}
[c]{cc}%
\subfigure[displacement at load\ increment~5 \label{figure:u55}] {\includegraphics[width=8cm]{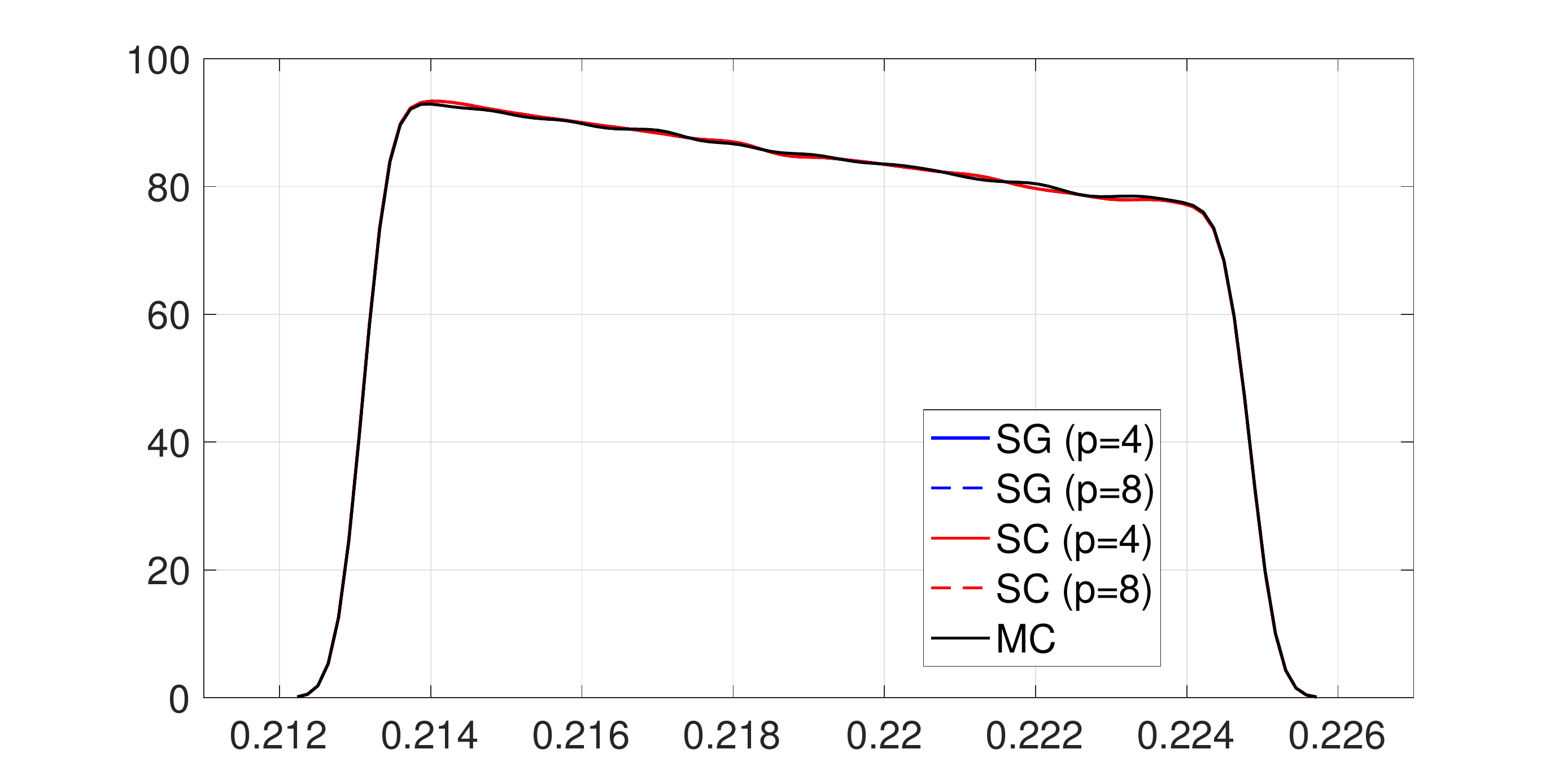}} &
\subfigure[displacement at load\ increment~10 \label{figure:u510}] {\includegraphics[width=8cm]{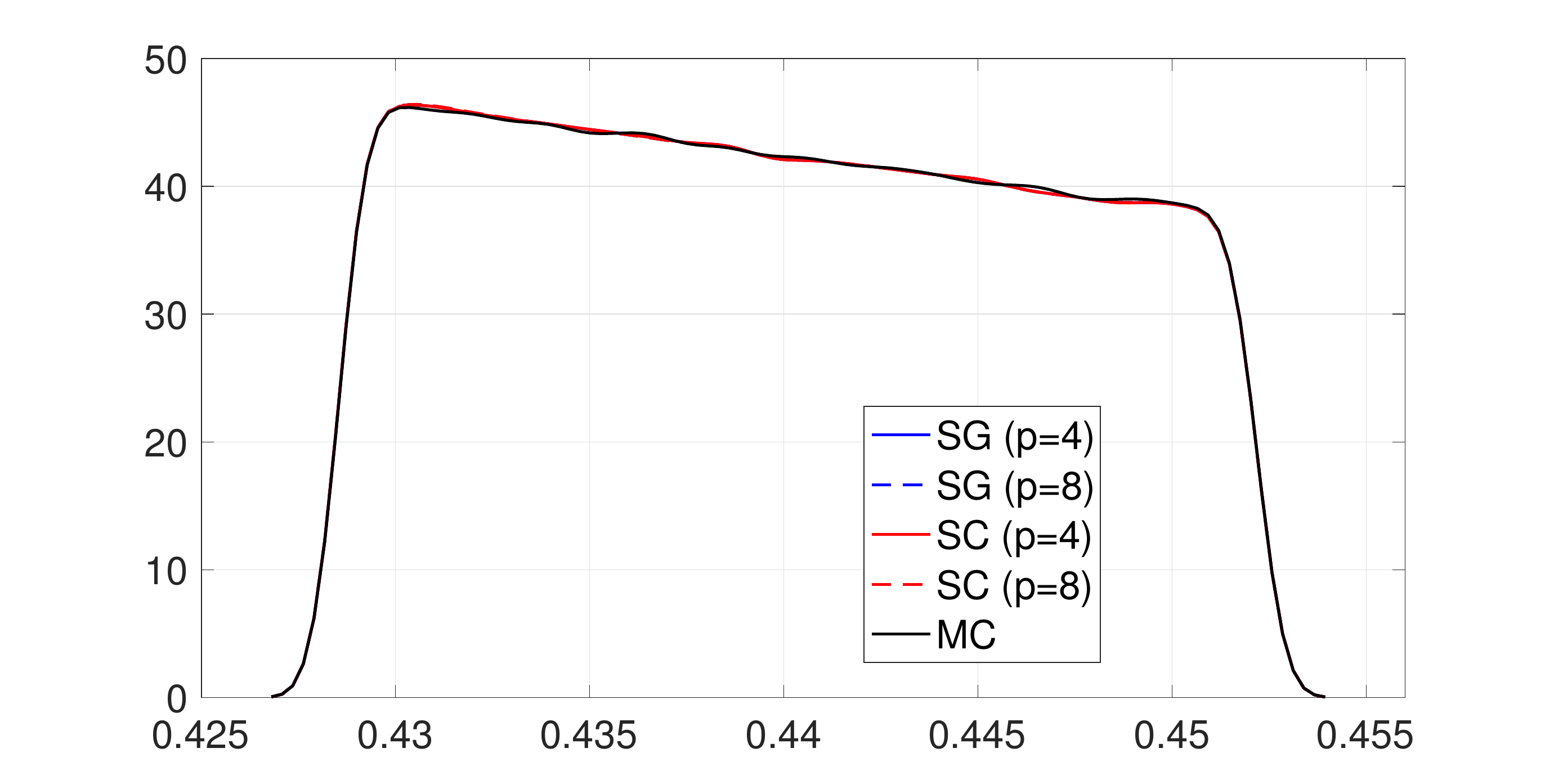}}\\
\subfigure[displacement at load\ increment~20 \label{figure:u520}] {\includegraphics[width=8cm]{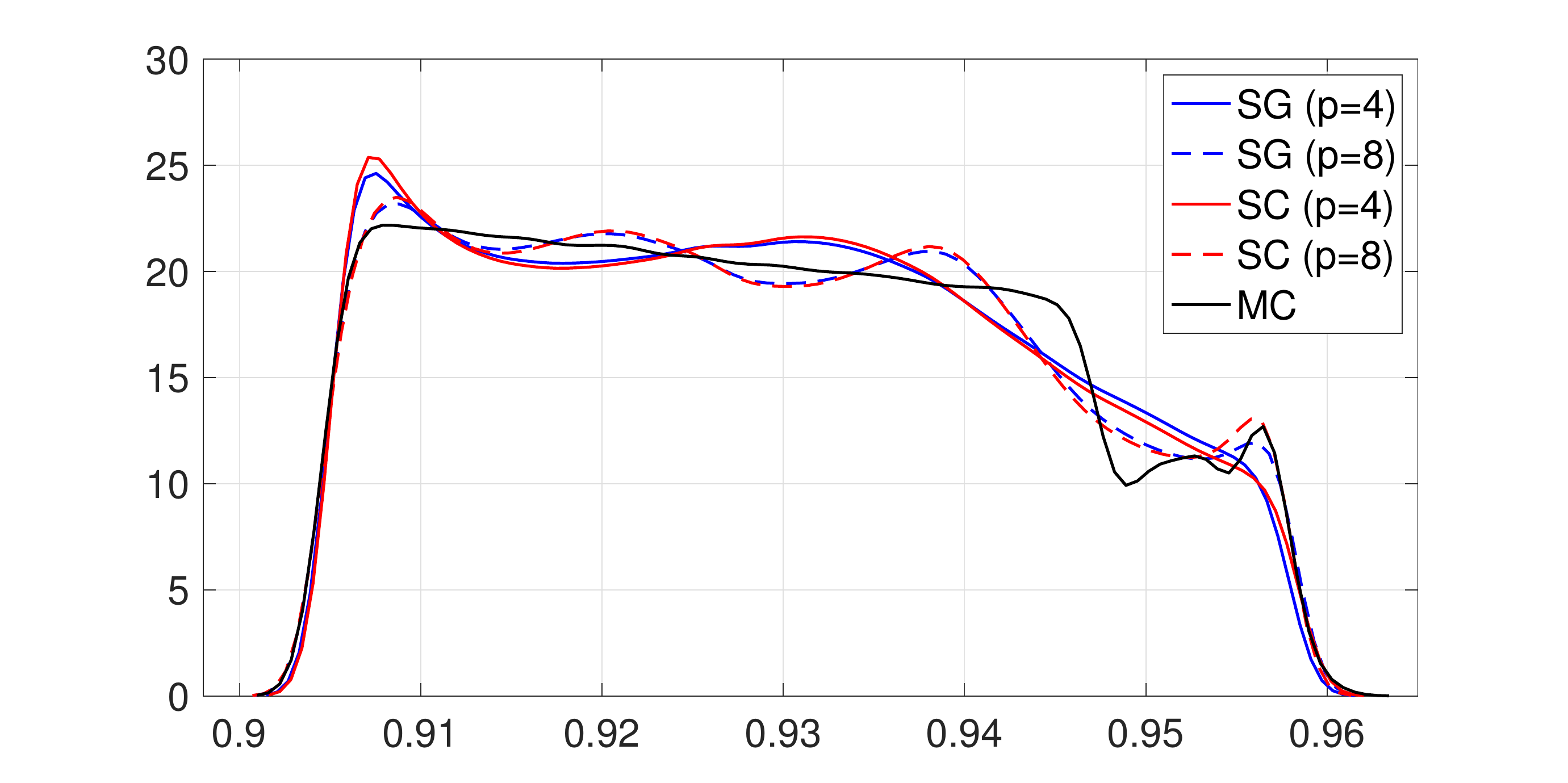}} &
\subfigure[displacement at load\ increment~50 \label{figure:u550}] {\includegraphics[width=8cm]{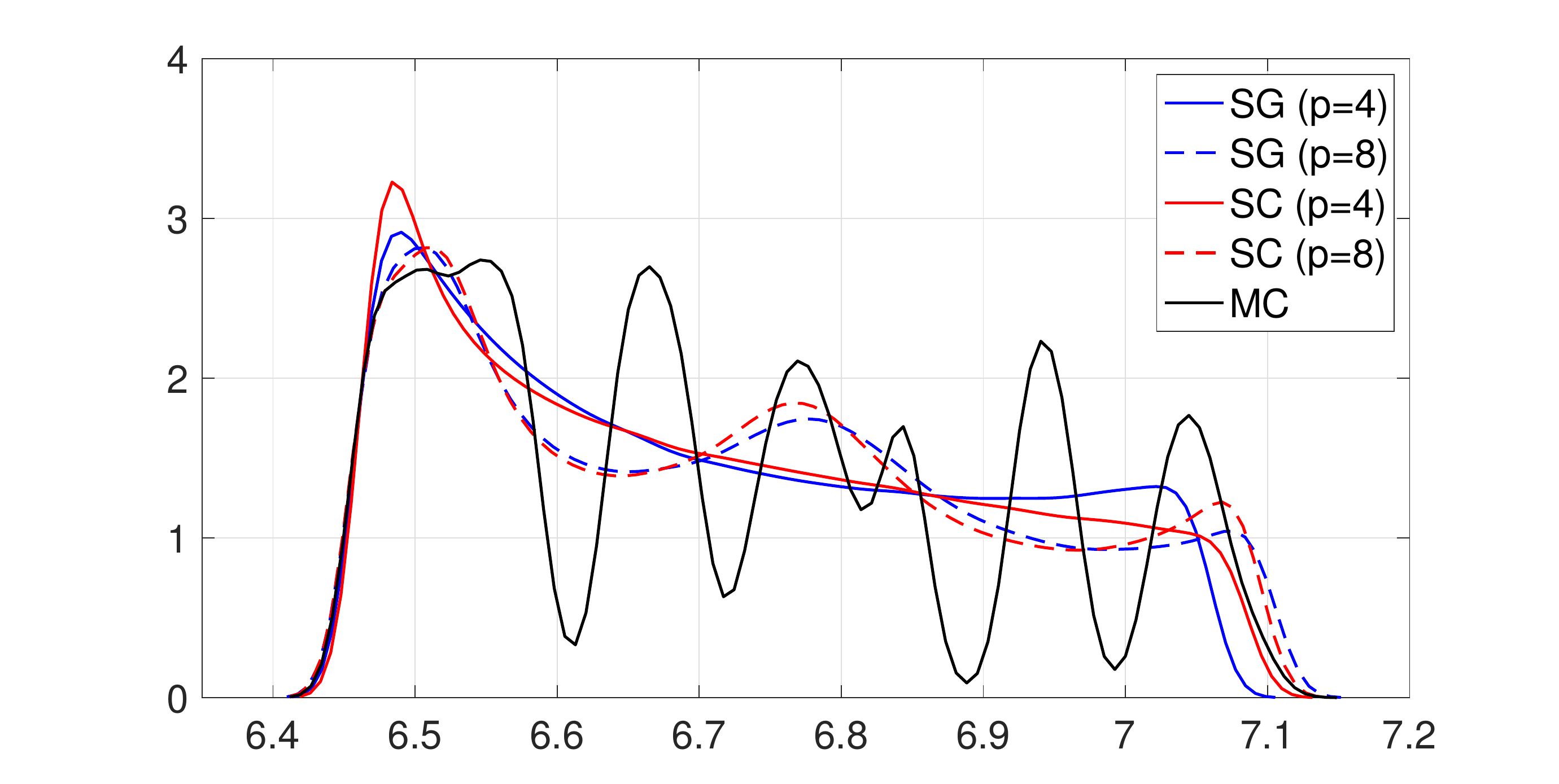}}\\
&
\end{tabular}
\caption{Estimated PDFs of the displacement at the location of the lateral
load at top of the shear wall during the loading.}
\label{figure:CoV5}%
\end{figure}

\begin{table}[b]
\caption{Displacement at the location of the lateral load at top of the shear wall: 
values of the mean~$\mu$ and standard deviation~$\sigma$ by Monte Carlo (MC), stochastic Galerkin (SG) and stochastic collocation (SC) methods with gPC degrees $p=4$ and $8$, comparison of the methods using root mean square error (RMSE), and estimated probability $\operatorname{Pr}(u-\widetilde{u}\geq0)$ in \%, where $u$ is the displacement from the stochastic problem and different setting of $\widetilde{u}$, in which $u_{m}$ is the displacement from the deterministic problem with the mean value parameters.}%
\label{tab:3}
\begin{center}%
\begin{tabular}
[c]{|c|c|c|c|c|c|c|}\hline
method & $\mu$ & $\sigma$ & $\operatorname{RMSE}$ & $\widetilde{u}=u_{m}$%
& $\widetilde{u}=1.01u_{m}$ & $\widetilde{u}=1.02u_{m}$\\\hline
\multicolumn{7}{|c|}{load increment 5}\\\hline
MC & $2.1876\times10^{-1}$ & $3.3858\times10^{-3}$ & - & $48.56 \,\%$ &
$30.76 \,\%$ & $13.27 \,\%$\\\hline
SG, $p=4$ & $2.1875\times10^{-1}$ & $3.3860\times10^{-3}$ & $4.7847\times
10^{-3}$ & $48.98\,\%$ & $30.70\,\%$ & $13.23\,\%$\\
SG, $p=8$ & $2.1875\times10^{-1}$ & $3.3860\times10^{-3}$ & $4.7847\times
10^{-3}$ & $48.98\,\%$ & $30.70\,\%$ & $13.23\,\%$\\\hline
SC, $p=4$ & $2.1875\times10^{-1}$ & $3.3860\times10^{-3}$ & $4.7847\times
10^{-3}$ & $48.98\,\%$ & $30.70\,\%$ & $13.23\,\%$\\
SC, $p=8$ & $2.1875\times10^{-1}$ & $3.3860\times10^{-3}$ & $4.7847\times
10^{-3}$ & $48.98\,\%$ & $30.70\,\%$ & $13.23\,\%$\\\hline
\multicolumn{7}{|c|}{load increment 10}\\\hline
MC & $4.3994\times10^{-1}$ & $6.8147\times10^{-3}$ & - & $49.09 \,\%$ &
$30.77 \,\%$ & $13.30 \,\%$\\\hline
SG, $p=4$ & $4.3993\times10^{-1}$ & $6.8147\times10^{-3}$ & $9.6302\times
10^{-3}$ & $48.99\,\%$ & $30.71\,\%$ & $13.25\,\%$\\
SG, $p=8$ & $4.3993\times10^{-1}$ & $6.8153\times10^{-3}$ & $9.6305\times
10^{-3}$ & $48.99\,\%$ & $30.71\,\%$ & $13.26\,\%$\\\hline
SC, $p=4$ & $4.3993\times10^{-1}$ & $6.8148\times10^{-3}$ & $9.6302\times
10^{-3}$ & $48.99\,\%$ & $30.71\,\%$ & $13.25\,\%$\\
SC, $p=8$ & $4.3993\times10^{-1}$ & $6.8153\times10^{-3}$ & $9.6305\times
10^{-3}$ & $48.99\,\%$ & $30.71\,\%$ & $13.26\,\%$\\\hline
\multicolumn{7}{|c|}{load increment 20}\\\hline
MC & $9.2863\times10^{-1}$ & $1.4732\times10^{-2}$ & - & $48.43 \,\%$ &
$29.91 \,\%$ & $12.25 \,\%$\\\hline
SG, $p=4$ & $9.2851\times10^{-1}$ & $1.4604\times10^{-2}$ & $2.0728\times
10^{-2}$ & $48.56\,\%$ & $29.07\,\%$ & $13.15\,\%$\\
SG, $p=8$ & $9.2861\times10^{-1}$ & $1.4755\times10^{-2}$ & $2.0835\times
10^{-2}$ & $48.31\,\%$ & $29.77\,\%$ & $13.25\,\%$\\\hline
SC, $p=4$ & $9.2853\times10^{-1}$ & $1.4648\times10^{-2}$ & $2.0759\times
10^{-2}$ & $48.65\,\%$ & $28.95\,\%$ & $13.18\,\%$\\
SC, $p=8$ & $9.2862\times10^{-1}$ & $1.4763\times10^{-2}$ & $2.0840\times
10^{-2}$ & $48.28\,\%$ & $29.78\,\%$ & $13.31\,\%$\\\hline
\multicolumn{7}{|c|}{load increment 50}\\\hline
MC & $6.6636\times10^{0}$ & $3.0779\times10^{-1}$ & - & $55.07 \,\%$ &
$45.23 \,\%$ & $31.54 \,\%$\\\hline
SG, $p=4$ & $6.7119\times10^{0}$ & $0.8756\times10^{-1}$ & $3.2090\times
10^{-1}$ & $53.02\,\%$ & $43.04\,\%$ & $33.91\,\%$\\
SG, $p=8$ & $6.7216\times10^{0}$ & $0.9085\times10^{-1}$ & $3.2220\times
10^{-1}$ & $56.25\,\%$ & $46.38\,\%$ & $34.91\,\%$\\\hline
SC, $p=4$ & $6.7130\times10^{0}$ & $0.8872\times10^{-1}$ & $3.2127\times
10^{-1}$ & $53.42\,\%$ & $43.17\,\%$ & $33.69\,\%$\\
SC, $p=8$ & $6.7212\times10^{0}$ & $0.9071\times10^{-1}$ & $3.2215\times
10^{-1}$ & $56.47\,\%$ & $46.41\,\%$ & $34.24\,\%$\\\hline
\end{tabular}
\end{center}
\end{table}

}}

\section{Iterative solution of the linearized systems}

\label{sec:precond} The solution of the linearized stochastic Galerkin systems
of equations~(\ref{equation.linearequequ2}) may be a computationally expensive
task due to its size, and also with respect to the large number of load
increments and steps, use of direct methods may be prohibitive. Therefore, we
also studied the problem of solving the stochastic Galerkin systems of
equations by iterative methods. Since the associated matrices are symmetric
and positive definite, we used conjugate gradient (CG)
method~\cite{babuvska2005solving}. However, the matrices are typically also
ill-conditioned, and construction of efficient preconditioners becomes an
important task for a practical implementation.
In this study, we use two preconditioners, the first one is the mean-based
preconditioner~\cite{pellissetti2000iterative, powell2009block} and the second
one is the hierarchical Gauss-Seidel
preconditioner~\cite{sousedik2014truncated}.

\subsection{Mean-based preconditioner}

\label{sec:precond:mean} The stiffness matrix derived from the SGFEM
formulation has a particular structure that can be exploited during the
process of solving the algebraic system of equations. In general, the
matrix$~K_{1}^{n}$ in equation~(\ref{equation.stiffdisc}) corresponds to the
mean properties of the system, and it has a much more considerable
contribution than the other $K_{i}^{n}$'s that represent random fluctuations
of the system from the mean, especially for smaller values of$~CoV$. Since the
submatrix $K_{1}^{n}$ contributes only to the block-diagonal of the stochastic
Galerkin matrix, the resulting system of linear equations exhibits a strong
block-diagonal dominance if the random properties represent only small
fluctuations from the mean values. Then, the mean-based preconditioner, given
in a matrix form as $I_{n_{\xi}}\otimes K_{1}^{n}$, will be a good
approximation to the stochastic Galerkin matrix$~\mathbf{K}^{n}$. This
block-diagonal matrix has the great advantage because it can be inverted by
inverting each block along the block-diagonal independently.
On the other hand, in systems with large random fluctuations, the off-diagonal
blocks have a much stronger contribution. For these cases the mean-based
preconditioner may not improve the convergence rate, because it does not
approximate the system matrix sufficiently.

\subsection{Hierarchical Gauss-Seidel preconditioner}

\label{sec:precond:hier} 
\begin{figure}[b]
\centering
{\includegraphics[width=10cm]{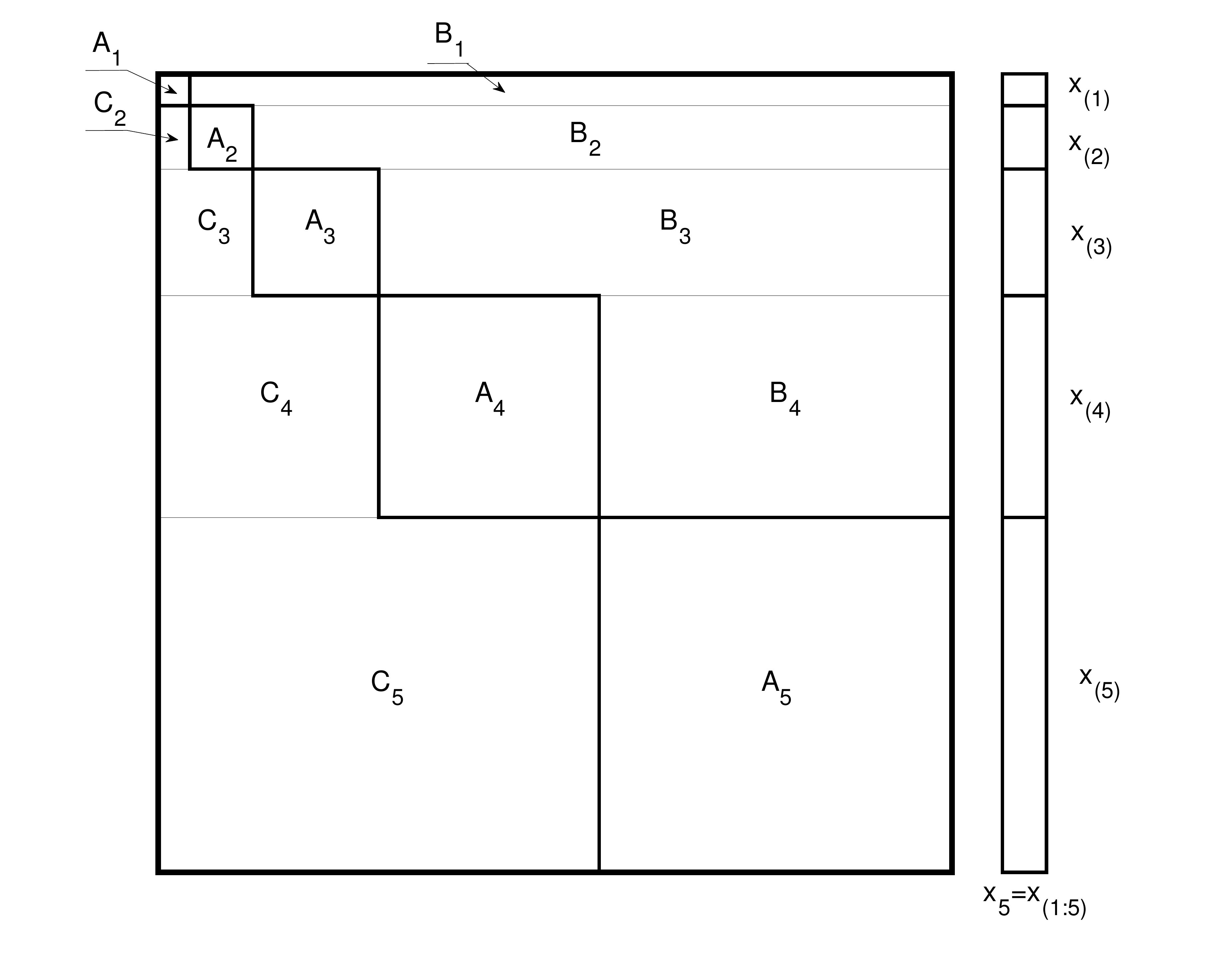}} \caption{The hierarchical
structure of the stochastic Galerkin matrix.}%
\label{fig:structure}%
\end{figure}

The block structure of the global stochastic
stiffness Galerkin matrix~(\ref{equation.stochstiff}) depends on the tensor
given by the values$~c_{jkm}$. We will now consider that the stochastic
Galerkin matrix$~\mathbf{K}^{n}$ in~(\ref{equation.stochstiff}) has a
hierarchical structure, cf. Fig.~\ref{fig:structure}, given as
\begin{equation}
\mathbf{K}^{n}=%
\begin{bmatrix}
A_{1} &  & B_{1} &  & \\
& \ddots &  &  & \\
C_{l} &  & A_{l} &  & B_{l}\\
&  &  & \ddots & \\
&  & C_{p+1} &  & A_{p+1}%
\end{bmatrix}
, \label{equation.hierstiff}%
\end{equation}
where $A_{1}=K_{1}^{n}$ is the matrix of the mean. The
decomposition~(\ref{equation.hierstiff}) is used to formulate the hierarchical
Gauss-Seidel preconditioner~\cite{sousedik2014truncated}.
A matrix-vector multiplication by the stochastic Galerkin matrix is in an
iterative solver performed using $K_{(k,m)}^{n}=\sum_{j=1}^{n_{\xi}}c_{jkm}%
{K}_{j}^{n}$, so we use only the constants $c_{jkm}$ and matrices ${K}_{j}%
^{n}$. The same strategy is used in the preconditioner for multiplications by
the submatrices $B_{l}$ and $C_{l}$, see~(\ref{equation.hierstiff}), and the
solves with the diagonal submatrices are approximated by the block-diagonal
solves with the mean matrix, that is $A_{l}\approx\widetilde{A}_{l}=I\otimes
A_{1}$ of the appropriate size. Finally, let us introduce the following
notation for a vector $x_{l}$, where $l=1,\dots,p+1$, cf.
Fig.~\ref{fig:structure}, as
\begin{equation}
x_{l}=x_{(1:l)}=%
\begin{bmatrix}
x_{(1)}\\
x_{(2)}\\
\vdots\\
x_{(l)}%
\end{bmatrix}
\label{equation.x}%
\end{equation}
The hierarchical Gauss-Seidel preconditioner (ahGS): $r_{p+1}\longmapsto
\upsilon_{p+1}$ for the system~(\ref{equation.linearequequ2}) is defined as
follows~\cite{sousedik2014truncated}:

Set the initial solution $v_{p+1}=0$ and update it in the following steps,
\begin{equation}
\upsilon_{(1)}=A_{1}^{-1}(r_{(1)}-B_{1}\upsilon_{(2:p+1)})
\label{equation.Hierarchical1}%
\end{equation}
\textbf{for} $l=2,\cdots,p$,
\begin{equation}
\upsilon_{(l)}=\widetilde{A}_{l}^{-1}(r_{(l)}-C_{l}\upsilon_{(1:l-1)}%
-B_{l}\upsilon_{(l+1:p+1)}) \label{equation.Hierarchical2}%
\end{equation}
\textbf{end}
\begin{equation}
\upsilon_{(p+1)}=\widetilde{A}_{p+1}^{-1}(r_{(p+1)}-C_{p}\upsilon_{(1:p)})
\label{equation.Hierarchical3}%
\end{equation}
\textbf{for} $l=p,\cdots,2$,
\begin{equation}
\upsilon_{(l)}=\widetilde{A}_{l}^{-1}(r_{(l)}-C_{l}\upsilon_{(1:l-1)}%
-B_{l}\upsilon_{(l+1:p+1)}) \label{equation.Hierarchical4}%
\end{equation}
\textbf{end}
\begin{equation}
\upsilon_{(1)}=A_{1}^{-1}(r_{(1)}-B_{1}\upsilon_{(2:p+1)})
\label{equation.Hierarchical5}%
\end{equation}
Since we initialize $v_{p}=0$, the multiplications by $B_{l}$, $l=1,\dots,p$,
vanish from~(\ref{equation.Hierarchical1}) and~(\ref{equation.Hierarchical2}%
).

\subsection{Numerical experiments}

\label{sec:precond:numerical}We tested the solvers using \textcolor{black}{both} the beam 
\textcolor{black}{and the shear wall}.
For solving the system of
equations~(\ref{equation.linearequequ2}), we used the CG method with the
mean-based (MB) and hierarchical Gauss-Seidel (ahGS)\ preconditioners from
Sections~\ref{sec:precond:mean} and~\ref{sec:precond:hier}, respectively. As
the stopping criterion we used a reduction of the 2-norm of the residual by a
factor$~10^{-8}$.

\begin{figure}[b]
\centering
\begin{tabular}
[c]{cc}%
{\includegraphics[width=8cm]{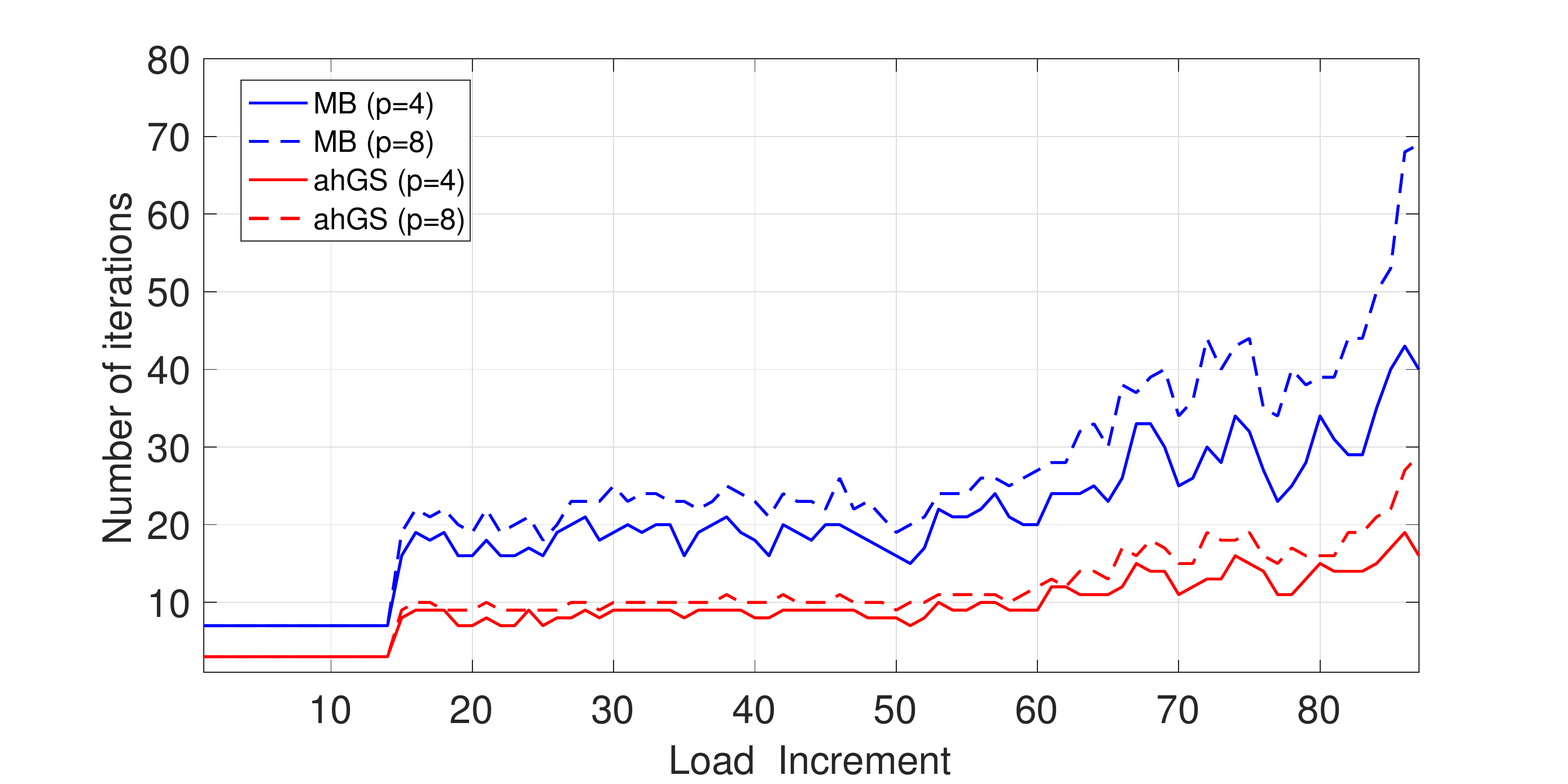}} &
{\includegraphics[width=8cm]{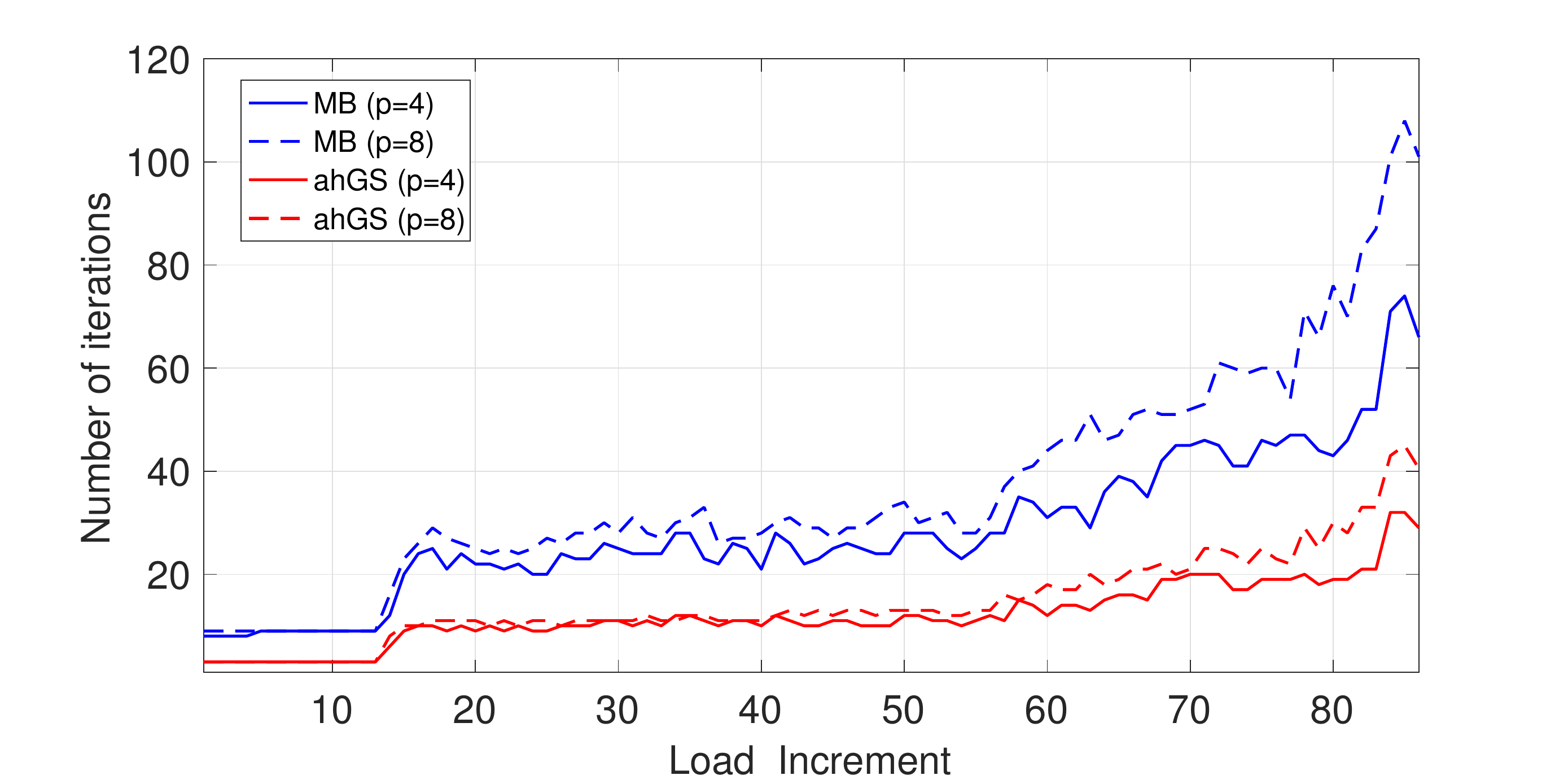}}\\
&
\end{tabular}
\caption{Iteration count for the beam at load increments with different preconditioners for
$CoV=5.77\%$ (left) and $10\%$ (right).}%
\label{figure:iteration2}%
\end{figure}

\begin{figure}[b]
\centering
\begin{tabular}
[c]{cc}%
\subfigure[$CoV\ 5.77\%$, increment~10 \label{figure:residualp6acc9cov577increment10}] {\includegraphics[width=4cm]{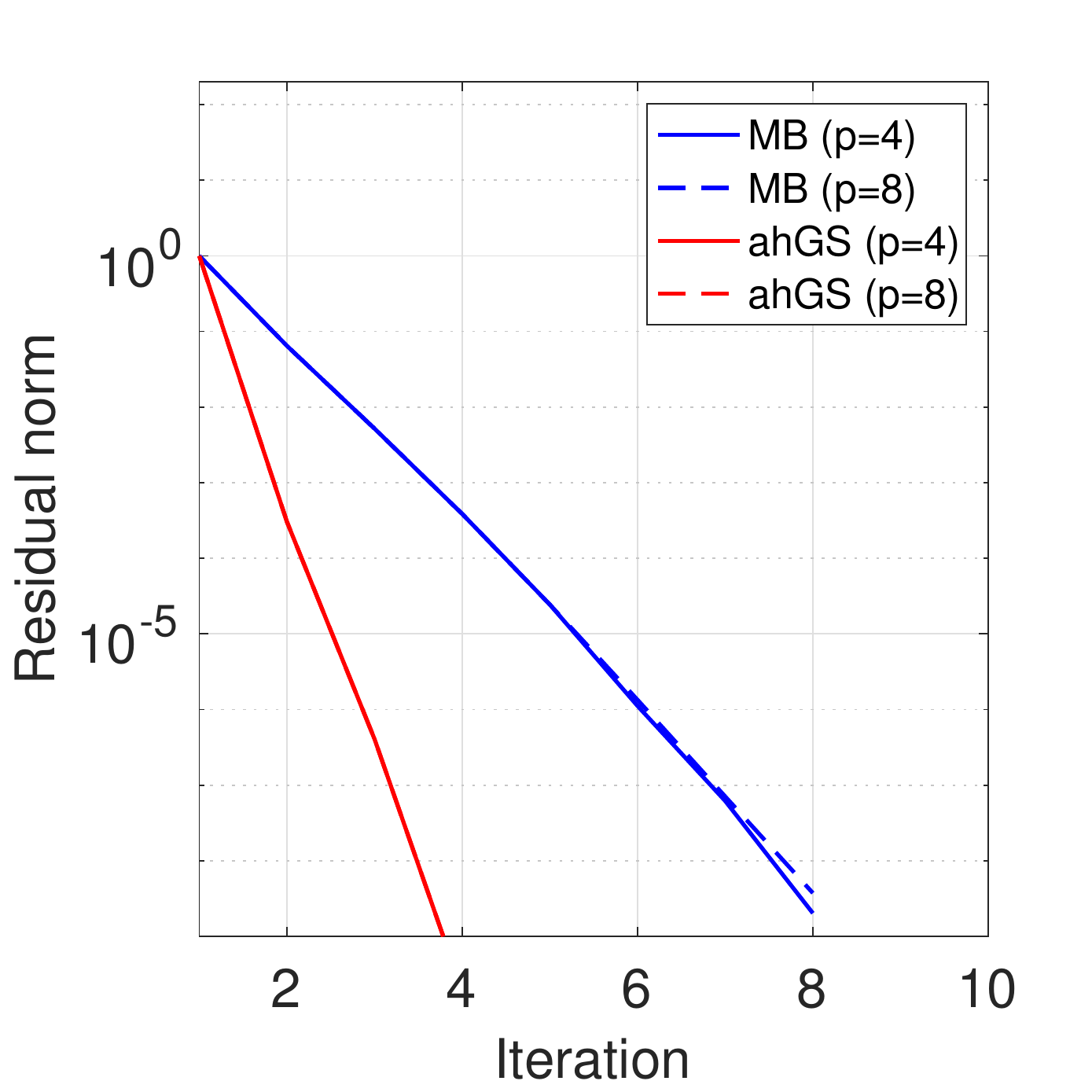}} &
\subfigure[$CoV\ 5.77\%$, increment~87 \label{figure:residualp6acc9cov577increment87}] {\includegraphics[width=8cm]{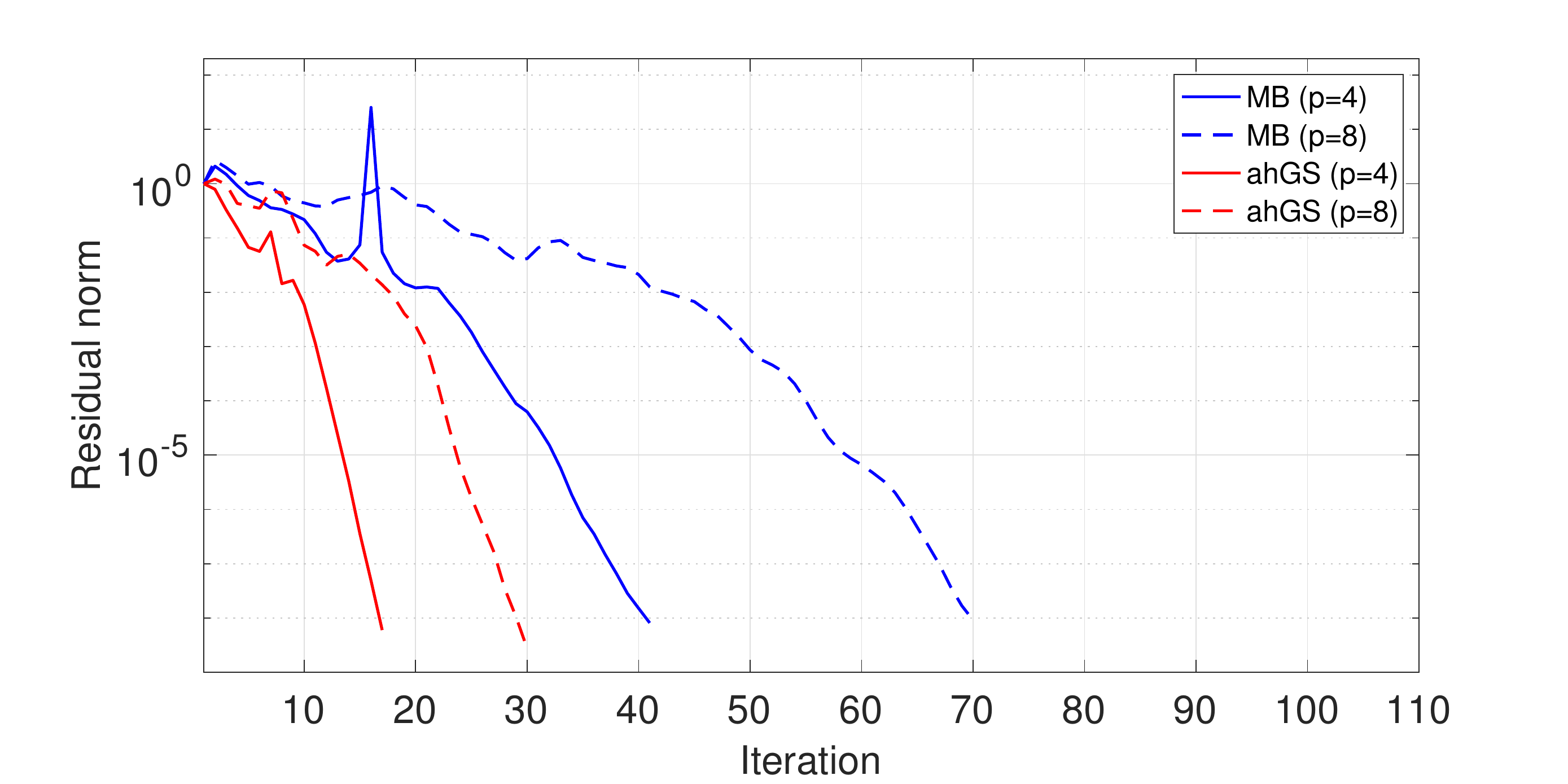}}\\
\subfigure[$CoV\ 10\%$, increment~10 \label{figure:residualp4acc10cov10increment10}] {\includegraphics[width=4cm,height=4cm]{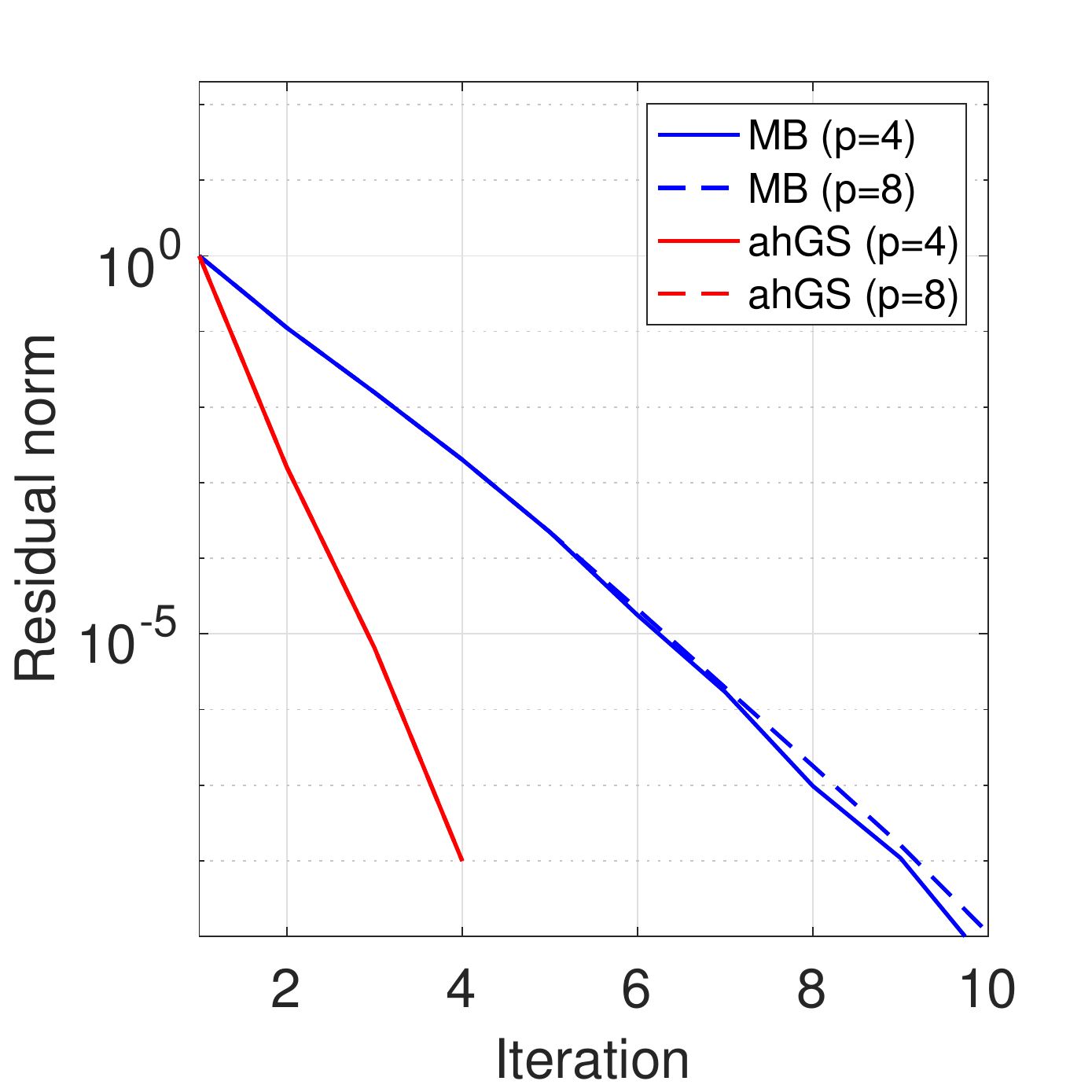}} &
\subfigure[$CoV\ 10\%$, increment~86 \label{figure:residualp4acc10cov10increment86}] {\includegraphics[width=8cm]{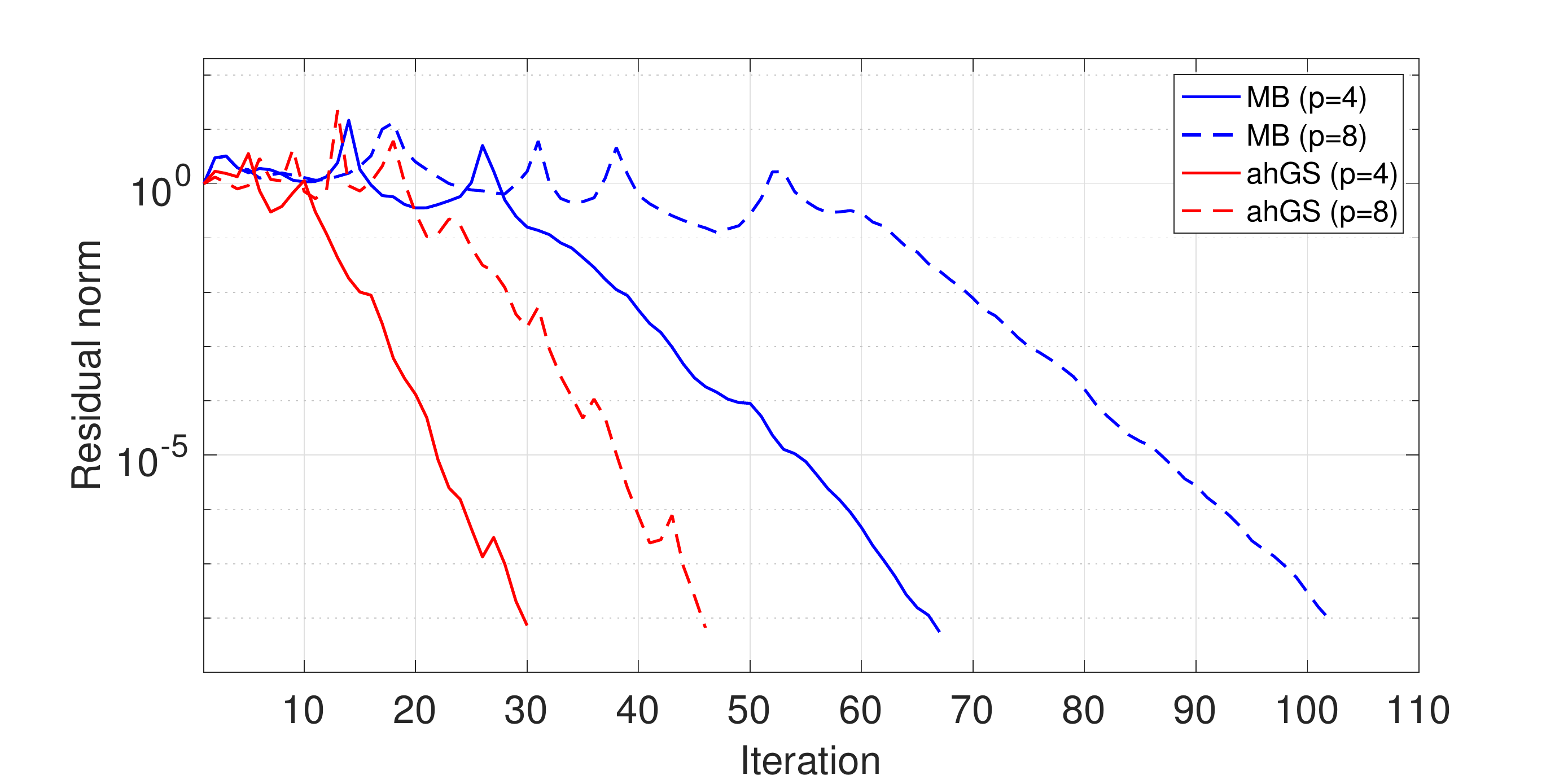}}\\
&
\end{tabular}
\caption{Residual norm at each iteration for the beam with $CoV=5.77\%$ (top) and $10\%$
(bottom).}%
\label{figure:residual}%
\end{figure}

\textcolor{black}{For the beam,} Figure~\ref{figure:iteration2} shows the numbers of iterations in each load
increment. We see that for all load increments the ahGS\ preconditioner
reduces the iteration count to approximately one half compared to the
MB\ preconditioner. For both preconditioners, the iteration count increases
slightly for the higher degree of the gPC\ polynomial, but the increase is
less pronounced for the ahGS preconditioner. The same observation can be made
for an increase of$~CoV$ by comparing the left and right panels in
Figure~\ref{figure:iteration2}. We can also easily deduce from
Figure~\ref{figure:iteration2} that the crack starts to develop in load
increment~$14$. In particular, the number of iterations is quite small and
constant during the initial load increments, but the character of the problem
changes and becomes challenging for the iterative solvers once the crack
starts to develop. The numbers of iterations slowly increase as the loading
progresses, but the overall growth of the iteration count is somewhat slower
for the ahGS\ preconditioner. Figure~\ref{figure:residual} then shows residual
history at load increment$~10$ and at the last load increments for the two
choices of $CoV$. This further illustrates the discussion above, and in
particular the efficiency gained due to the use of the ahGS\ preconditioner
throughout the entire loading process.

\begin{figure}[b]
\centering
\includegraphics[width=8cm]{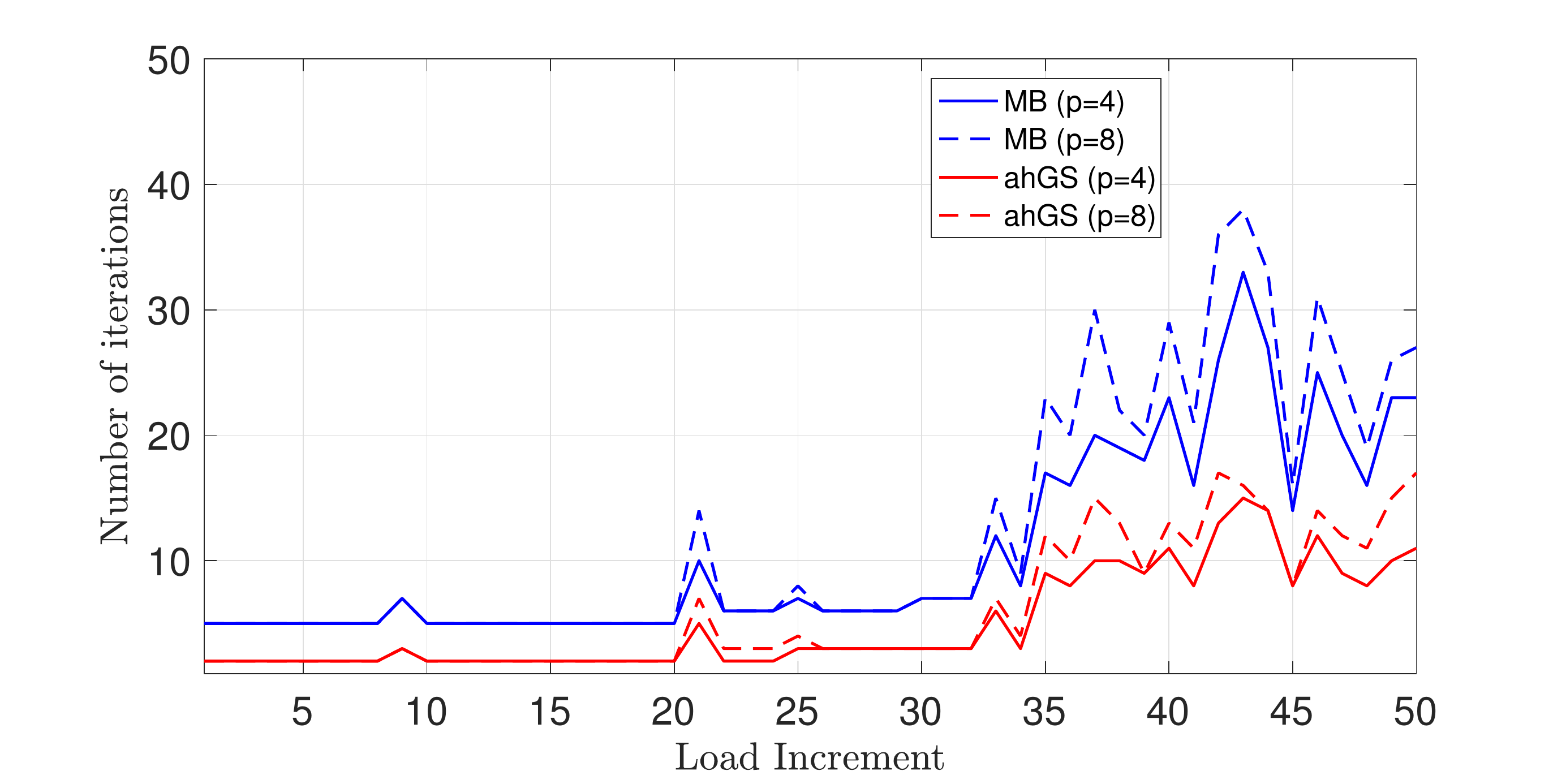} 
\caption{Iteration count at load increments with different preconditioners for the shear wall.}%
\label{figure:iter2}%
\end{figure}

\begin{figure}[b]
\centering
\begin{tabular}
[c]{cc}%
\subfigure[increment~5 \label{figure:residualcov5inc5}] {\includegraphics[width=8cm]{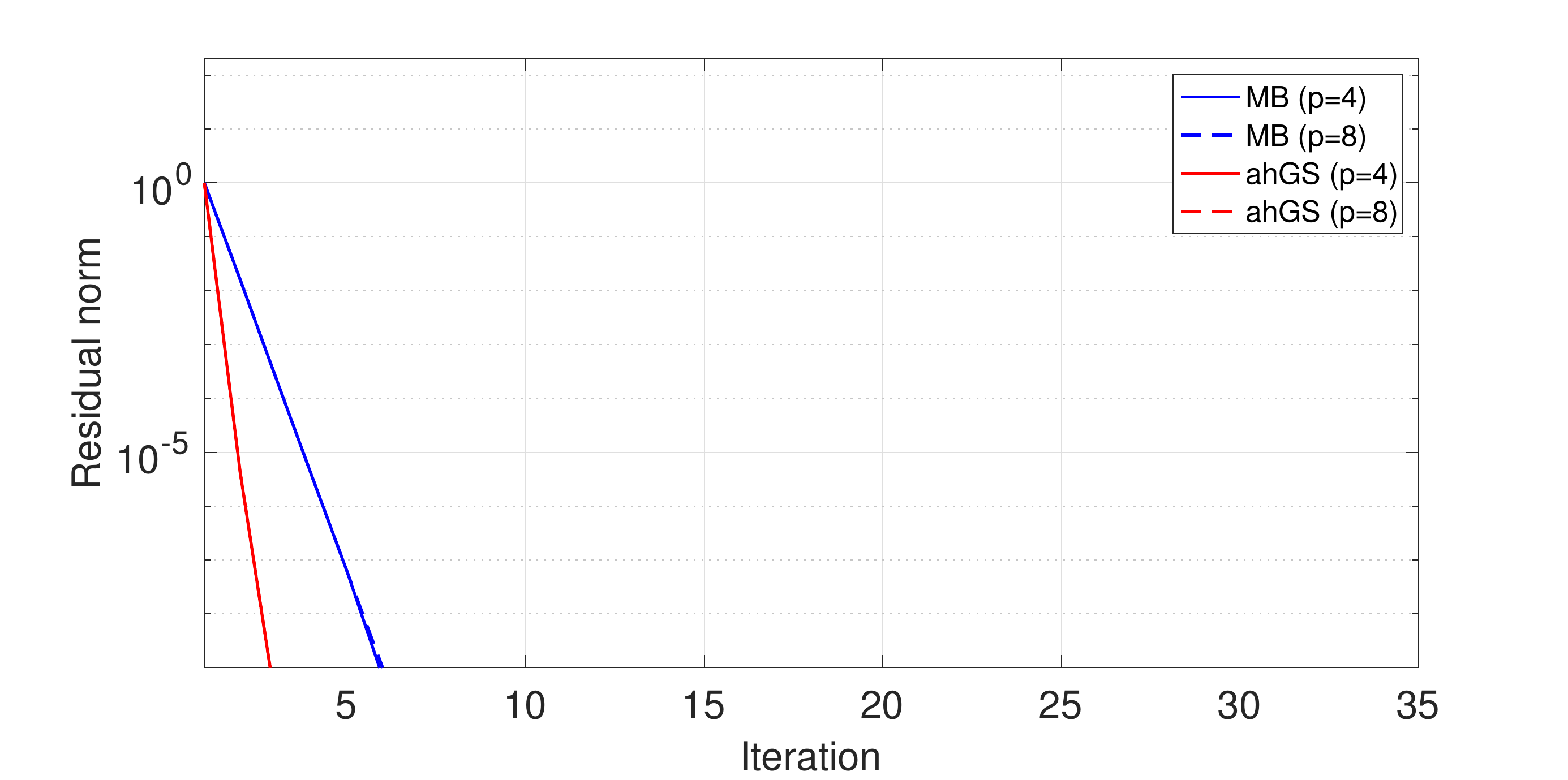}} & 
\subfigure[increment~10 \label{figure:residualcov5inc10}] {\includegraphics[width=8cm]{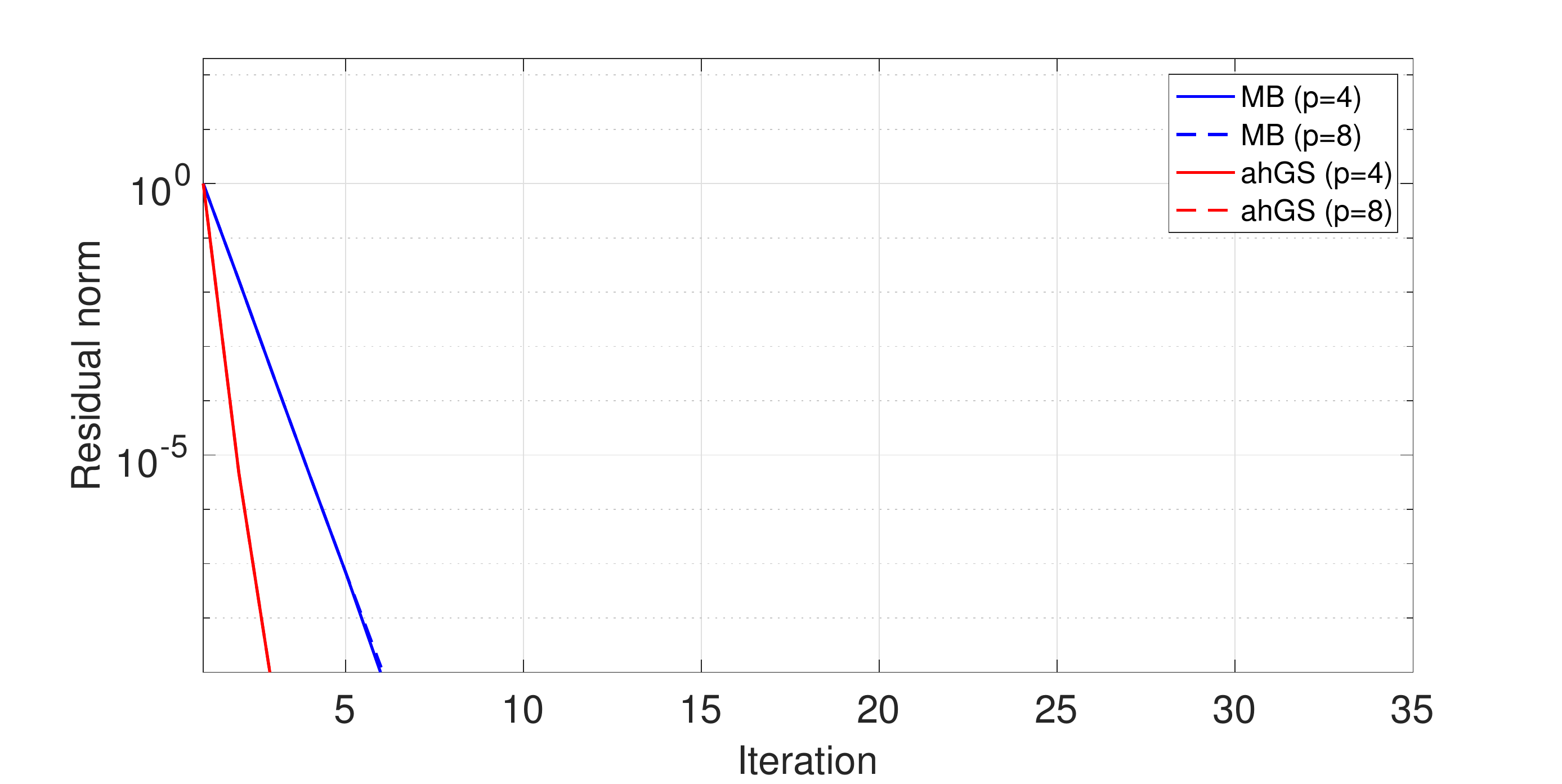}}\\
\subfigure[increment~20 \label{figure:residualcov5inc20}] {\includegraphics[width=8cm]{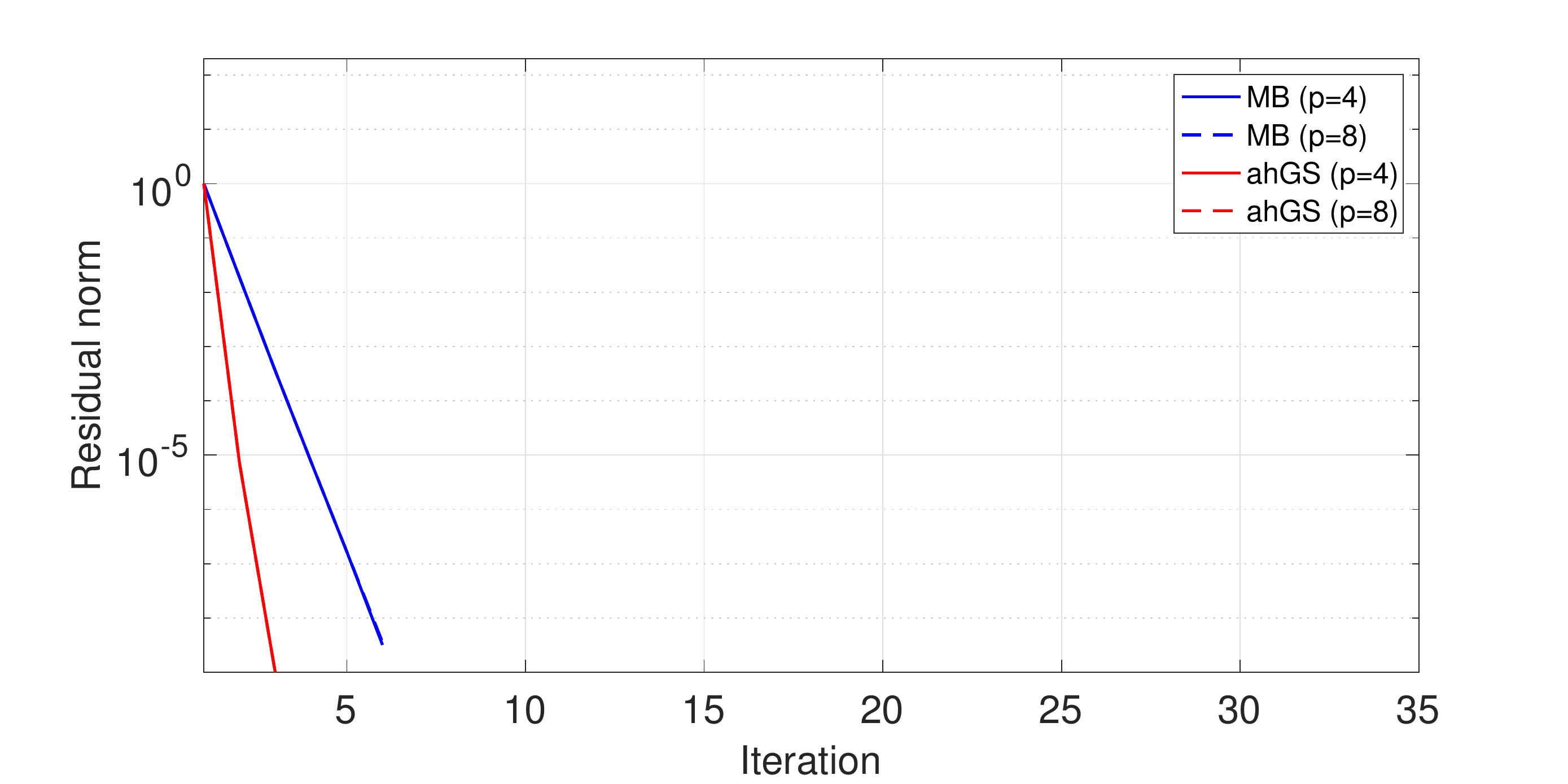}} &
\subfigure[increment~50 \label{figure:residualcov5inc50}] {\includegraphics[width=8cm]{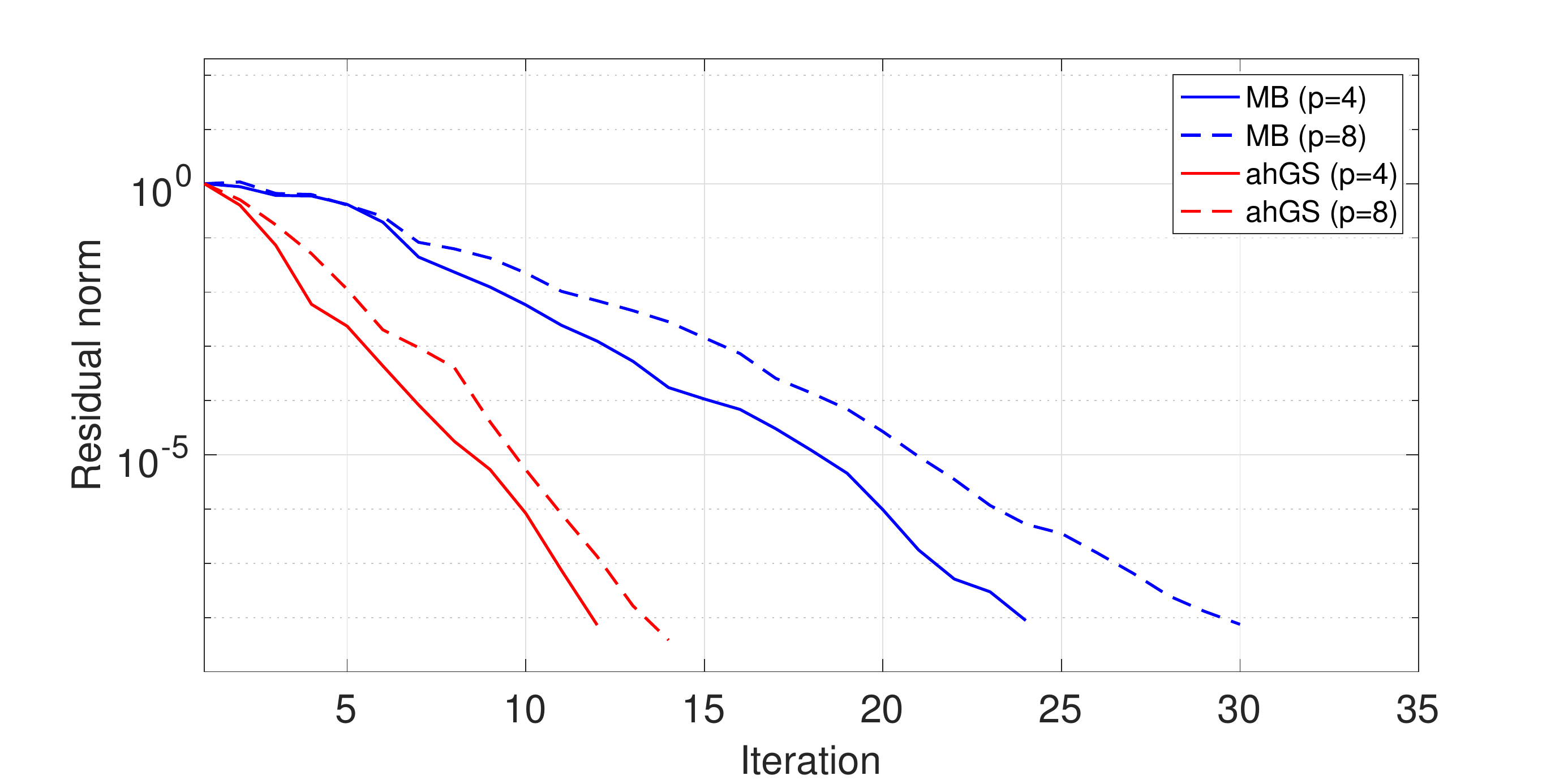}}\\
&
\end{tabular}
\caption{Residual norm at each iteration for the shear wall.}%
\label{figure:residua2}%
\end{figure}

{\color{black}{%
For the shear wall, the 
Figure~\ref{figure:iter2} shows the numbers of iterations in each load increment. We see, similarly as for the beam, that for all load increments the ahGS\ preconditioner reduces the iteration count compared to the MB\ preconditioner to approximately one half. 
And in a similar observation for both preconditioners, the iteration count increases slightly for the higher degree of the gPC\ polynomial, but the increase is less pronounced for the ahGS preconditioner.   
Figure~\ref{figure:residua2} then shows residual history at load increment$~5$, $~10$, $~20$ and~$50$, which is the last load increment. 
The plots further illustrate the discussion above, and in particular the efficiency gained due to the use of the ahGS\ preconditioner throughout the entire loading process.
}}

\section{Conclusion}

\label{sec:conclusion} A new methodology for the extension of the stochastic
Galerkin finite element method (SGFEM) is presented to efficiently propagate
the uncertainty in nonlinear elasticity problems. It can be applied, in
general, to models of reinforced concrete members with random material
properties. In the nonlinear SGFEM, the linearization scheme based on the
modified Newton-Raphson method. By directly updating the gPC expansions of the
strain and stress vectors between each loading step and updating the stiffness
matrix at the first step of load increments, the stochastic displacements at
each level of loading are obtained. The performance of the nonlinear SGFEM is
tested using a reinforced concrete beam with random initial concrete modulus of elasticity \textcolor{black}{and a shear wall with random 
maximum compressive stress of concrete.}
We illustrate by numerical experiments that even a low-degree
gPC\ polynomial provides a smooth interpolant to the full solution provided by
Monte Carlo simulation, and it provides an accurate estimate of several
statistical indicators associated with the probability distribution of the
response of the structure.
Since the linear systems associated with the use of the stochastic Galerkin
method may be very large and many of them need to be solved due to the
incremental loading, we also studied their iterative solution using
preconditioned conjugate gradient method with mean-based (MB) and hierarchical
Gauss-Seidel (ahGS) preconditioners. Numerical experiments show that the ahGS
preconditioner reduces the number of iterations to approximately one half of
iterations needed using the MB preconditioner.
Using the proposed methodology, it becomes practical to accurately approximate
the distribution of the structural response, which could be used for risk
assessment and more efficient engineering design of nonlinear structures that
include uncertainty in the material properties.


\acknowledgements 
We would like to thank the anonymous referees for their comments and suggestions.
B. Soused\'{\i}k was supported by the U.S. National Science
Foundation under award \mbox{DMS1913201}. Part of the work was completed while
M.~S. Ghavami was visiting University of Maryland, Baltimore County.

\bibliographystyle{IJ4UQ_Bibliography_Style}
\bibliography{paper_msgbshd}

\end{document}